\documentclass[journal]{IEEEtran}

\usepackage{amsmath}
\usepackage{amssymb}
\usepackage{mathrsfs}
\usepackage{turnstile}
\usepackage{graphicx}
\usepackage{color}
\usepackage{tensor}
\usepackage{algorithmic}
\usepackage{cite}
\usepackage{subfigure}

\usepackage[caption=false,font=footnotesize]{subfig} 

\providecommand {\pg} [1] { \left\{  #1 \right \} }
\providecommand {\pt} [1] { \left(  #1 \right) }
\providecommand {\pq} [1] { \left[  #1 \right] }

\providecommand{\fr} [2] {\frac{#1}{#2} \,}

\newcommand {\re} [1] {\mathrm{Re} \pg{#1} }
\newcommand {\im} [1] { \mathrm{Im} \pg{#1} }

\newcommand {\jrm}  {\mathrm{j}}
\newcommand {\trm}  {\mathrm{t}}
\newcommand {\esp} [1] {\mathrm{e}^{#1}}

\newcommand {\uv}[1] {\mathbf{u}_{#1}}
\newcommand {\g}[1] {\mathbf{#1}}

\newcommand {\om} {\omega}

\newcommand {\ii} {\infty}
\newcommand {\dd} {\, \textrm{d}}

\newcommand {\kx} {k_{x}}
\newcommand {\ky} {k_{y}}
\newcommand {\kz} {k_{z}}

\newcommand {\ux} {\mathbf{u}_{x}}
\newcommand {\uy} {\mathbf{u}_{y}}
\newcommand {\uz} {\mathbf{u}_{z}}

\newcommand {\intii} {\int_{-\ii}^{+\ii}}

\begin{document}

\title{Spectral-Domain Method of Moments\\ for the Modal Analysis of Line Waveguides}
%
%
%

\author{Giampiero~Lovat,
      Walter~Fuscaldo,~\IEEEmembership{Senior~Member,~IEEE,}
      Massimo~Moccia,\\
      Giuseppe~Castaldi,
      Vincenzo~Galdi,~\IEEEmembership{Fellow,~IEEE,}
      and~Paolo~Burghignoli,~\IEEEmembership{Senior~Member,~IEEE}
\thanks{G. Lovat and P. Burghignoli are with the DIAEE and DIET Departments, respectively, of Sapienza University of Rome, 00184 Rome, Italy (e-mail: giampiero.lovat@uniroma1.it, paolo.burghignoli@uniroma1.it).}
\thanks{W. Fuscaldo is with the Institute for Microelectronics and Microsystems, National Research Council (CNR-IMM), 00133 Rome, Italy (email: walter.fuscaldo@cnr.it).}
\thanks{M. Moccia, G. Castaldi and V. Galdi are with the Fields \& Waves Laboratory, Department of Engineering, University of Sannio, I-82100 Benevento, Italy (email: massimo.moccia@unisannio.it, castaldi@unisannio.it, vgaldi@unisannio.it).}
}

%
%

\markboth{}%
{Lovat \MakeLowercase{\textit{et al.}}: Spectral-Domain Method of Moments for Line Waveguides}

\maketitle

\begin{abstract}
A rigorous full-wave modal analysis based on the method of moments in the spectral domain is presented for line waveguides constituted by two-part impedance planes with arbitrary anisotropic surface impedances. An integral equation is formulated by introducing an auxiliary current sheet on one of the two half planes and extending the impedance boundary condition of the complementary half plane to hold on the entire plane. The equation is then discretized with the method of moments in the spectral domain, by employing exponentially weighted Laguerre polynomials as entire-domain basis functions and performing a Galerkin testing. Numerical results for both \emph{bound} and \emph{leaky} line waves are presented and validated against independent results, obtained for isotropic surface impedances with the analytical Sommerfeld--Maliuzhinets method and for the general anisotropic case with a commercial electromagnetic simulator. The proposed approach is computationally efficient, can accommodate the presence of spatial dispersion, and offers physical insight into the modal propagation regimes.
\end{abstract}


\IEEEpeerreviewmaketitle

\section{Introduction}

\IEEEPARstart{C}{ontrolling} 
the propagation of surface waves (SWs) has been a significant area of focus in microwave and antenna engineering for a long time \cite{Yang:2019se} and is now gaining renewed attention due to the emergence of low-dimensional materials, both artificial (metasurfaces \cite{Holloway:2012ao}) and natural (e.g., graphene) \cite{Bao:20192D}. These platforms allow for precise control of the surface-impedance properties, thereby enabling advanced SW-manipulation strategies and new effects, including waveguiding \cite{Gregoire:2011sw,Quarfoth:2013at}, transformation electromagnetics (EM) \cite{Vakil:2011to,Patel:2014te,Mencagli:2014mb}, and hyperbolic dispersion \cite{Yermakov:2015hw,Gomez-Diaz:2016fo}.

Within this framework, even simple planar discontinuities in the surface impedance can exhibit interesting properties. For example, they can support edge modes that are localized both in-plane and out-of-plane, guiding the EM power along a one-dimensional (1-D) path. These so-called ``line waves'' (LWs) represent a reduced-dimensional form of conventional SWs, and have been demonstrated to occur in planar junctions of metasurfaces with dual capacitive/inductive properties \cite{Horsley:2014od,dia2017guiding} or {\em parity-time} (${\cal PT}$) symmetry (i.e., balanced gain/loss) \cite{moccia2020line}.

Line waves exhibit a range of highly desirable characteristics, including 
strong localization and field enhancement, broad bandwidth, polarization-dependent propagation, topological robustness, spectral degeneracies, and bound/leaky regimes \cite{dia2017guiding,Bisharat:2019ed,moccia2020line,Moccia:2021ep,Xu:2021lw,Moccia:2023lw}. These features make them a promising option for various applications in the fields of antennas, integrated photonics, and optical sensing. Furthermore, the possibility of using gate-tunable graphene platforms at THz frequencies has been proposed \cite{Bisharat:2018ml}, which offers the potential for dynamic reconfiguration of wave pathways, confinement, and polarization.

Compared to SWs, the modeling of LWs is remarkably more challenging. In the simplest (isotropic, single-junction) scenarios, exact analytical solutions can be found by relying on the Sommerfeld--Maliuzhinets  \cite{kong2019analytic} or Wiener--Hopf \cite{Kay:1959so} techniques
commonly utilized in diffraction theory, but they are complicated and lack a physically incisive parameterization. 
General numerical methods, such as finite elements \cite{COMSOL:2015}, can be used for comprehensive parametric studies of more complex configurations, such as anisotropic or multi-junction systems. However, these methods still lack a clear physical understanding and do not easily accommodate important aspects such as improper (complex-valued) eigenmodes and spatial-dispersion effects. As a result, a comprehensive understanding and categorization of surface-impedance properties that support LWs in bound or leaky regimes is still lacking.

Against the above-mentioned limitations, this paper presents a rigorous full-wave modal analysis of LWs (preliminary results have been presented in \cite{lovatAPS23}). The proposed approach is based on the spectral-domain method of moments (MoM), and is capable of accommodating anisotropy and spatial dispersion. Furthermore, it can naturally handle both proper and improper eigenmodes. Our proposed method provides a computationally effective and physically insightful approach for the study of LWs.

The paper is structured as follows. 
Section \ref{Sec:IEF} outlines the problem formulation in terms of an integral equation, and Sec. \ref{Sec:MoM} details its MoM-based numerical solution.
In Sec. \ref{sec:results}, representative numerical results are discussed and compared with theoretical \cite{kong2019analytic} and independent numerical \cite{COMSOL:2015} predictions to validate and calibrate the proposed approach.
Finally, in Section \ref{Sec:Conclusions}, brief concluding remarks and suggestions for future research are provided. Additional theoretical and simulation details can be found in the Appendices.

\begin{figure}[!t]
	\centering
	\includegraphics[width=\columnwidth]{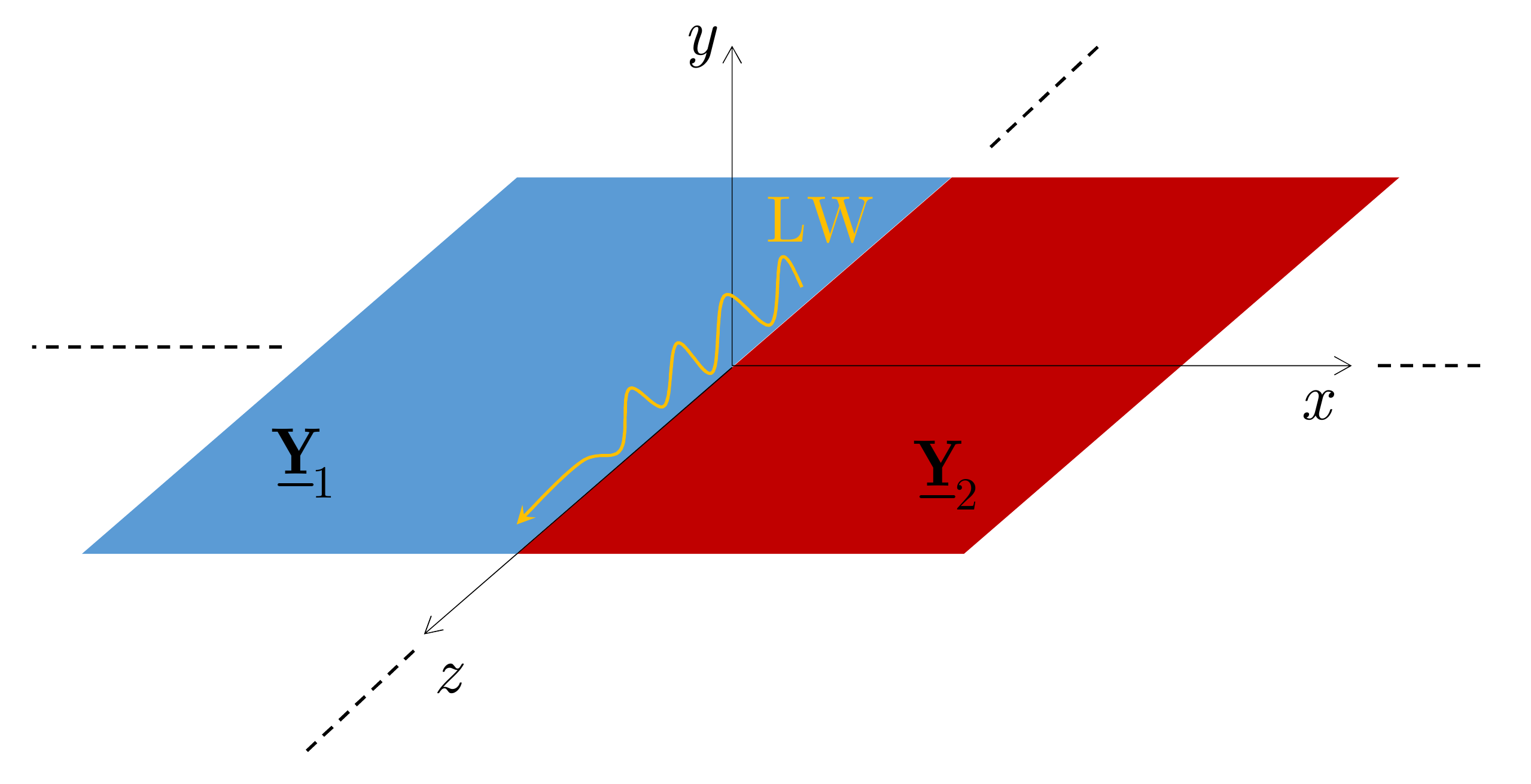}
	\caption{A two-part impedance plane, constituted by two half planes with dyadic surface admittances $\underline{\g{Y}}_{1,2}$, supporting a line wave (LW), with the relevant Cartesian reference system.} \label{fig:struct} 
\end{figure}

\section{Integral-Equation Formulation}
\label{Sec:IEF}

\subsection{Structure Description and Line-Wave Modal Field}

The configuration under analysis is the two-part impedance plane shown in Fig. \ref{fig:struct} and characterized by the surface impedance conditions
\begin{subequations}
\begin{align} 
   \uv{y} \times \g{H} &= \underline{\g{Y}}_1 \cdot \g{E}_\tau \,,\quad y=0, \; x<0 \,,\\
   \uv{y} \times \g{H} &= \underline{\g{Y}}_2 \cdot \g{E}_\tau \,,\quad y=0, \; x>0 \,,
\end{align}\label{eq:sibc}\end{subequations}%
where $\uv{y}$ is the unit vector of the $y$-axis, $\underline{\g{Y}}_{1,2}$ are the surface admittance dyadics of the two half planes $x<0$ and $x>0$, respectively, and the subscript $\tau$ indicates the component of the electric field tangential to the plane $y=0$. Note that the plane $y=0$ is an opaque boundary for the EM field, which is nonzero only in the vacuum half space $y>0$.

We aim at finding modal fields\footnote{As is known, the term `mode' is ambiguous, as it can be used both in `eigenwave' problems (in which the eigenvalue is the complex wavenumber of a traveling wave oscillating with an assigned real-valued radian frequency) and `eigenfrequency' problems (in which the eigenvalue is a complex-valued radian frequency). Throghout the paper, the terms `mode' and `modal' refer to eigenwaves supported by the two-part impedance plane.} in a time-harmonic regime $\esp{\jrm \om t}$ where $\omega$ is an assigned real-valued radian frequency, i.e., solutions of the homogeneous Maxwell's equations which satisfy the boundary conditions \eqref{eq:sibc}, having the form
\begin{align}
\label{eq:emod}
   \g{E}\pt{x,y,z}  &= \g{e}\pt{x,y} \esp{-\jrm \kz z} \,, \\
   \label{eq:hmod}
   \g{H}\pt{x,y,z}  &= \g{h}\pt{x,y} \esp{-\jrm \kz z} \,,
\end{align}
where $\kz = \beta_z - \jrm \alpha_z$ is the (generally complex) longitudinal propagation wavenumber of the LW, i.e., the eigenvalue of our modal problem, whereas $\g{e},\g{h}$ are the modal vectors, i.e., the relevant eigenvectors.

%
%
%


\subsection{Integral-Equation Formulation}

\begin{figure}[!t]
	\centering
	\includegraphics[width=\columnwidth]{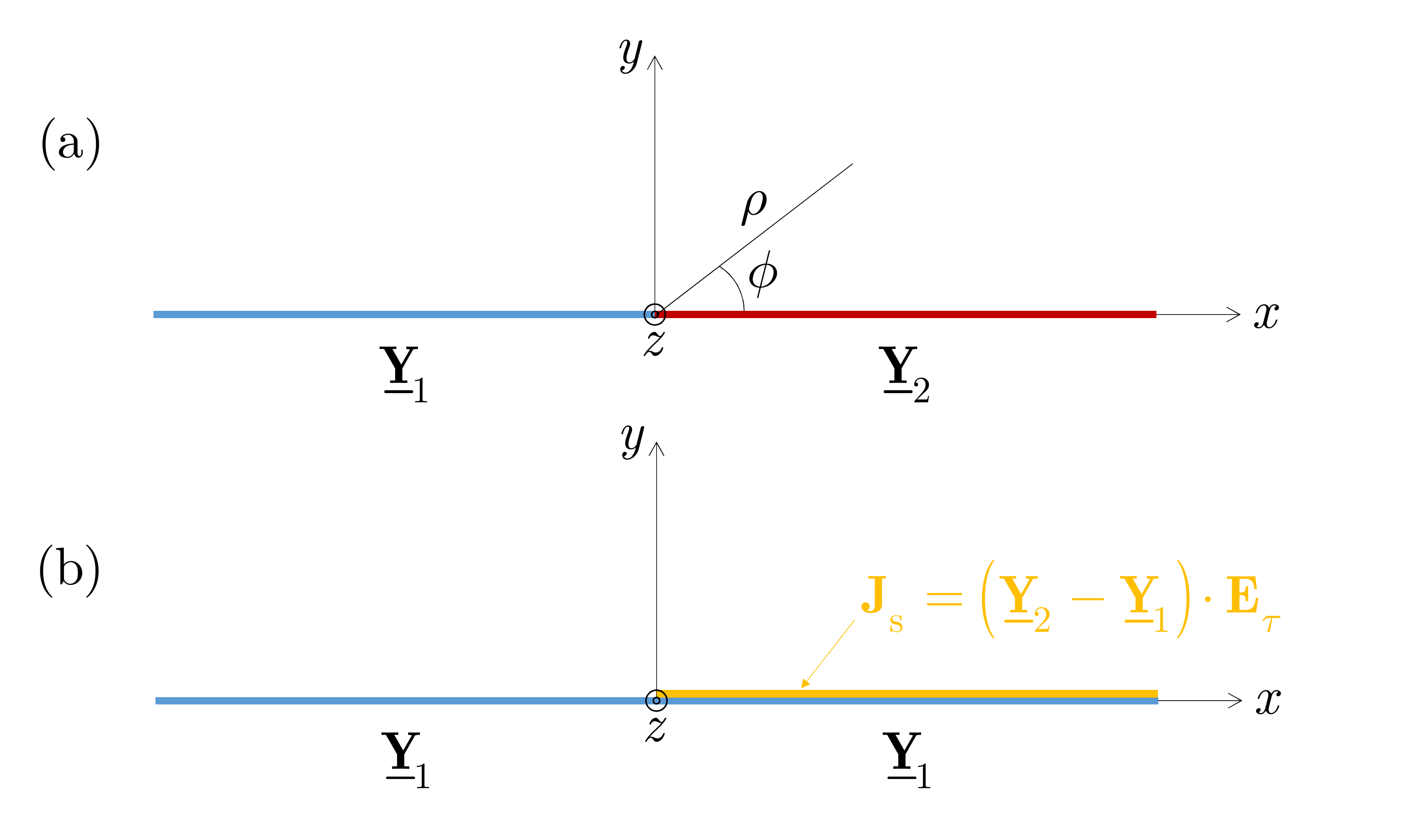}
	\caption{(a) Original two-part impedance configuration. (b) Equivalent configuration made of a uniform impedance plane and an electric current sheet.} \label{fig:equiv} 
\end{figure}

The modal problem will be formulated through an integral equation by replacing the original configuration of Fig. \ref{fig:equiv}(a) with the equivalent configuration shown in Fig. \ref{fig:equiv}(b), where the impedance boundary condition involving the  surface admittance dyadic $\underline{\g{Y}}_1$ has been extended to hold on the entire plane $y=0$:
\begin{align} \label{eq:sibc2}
   \uv{y} \times \g{H} &= \underline{\g{Y}}_1 \cdot \g{E}_\tau \,,\quad y=0,
\end{align}
and an auxiliary electric-current sheet with current density $\g{J}_{\mathrm{s}}$ has been introduced on the half plane $y=0,x>0$. Assuming that this behaves as a \textit{resistive sheet} \cite{bleszynski1993surface}, characterized by the transition conditions %
\begin{align} \label{eq:trans1}
   &\g{E}_\tau^+ = 
 \g{E}_\tau^- =: \g{E}_\tau \,,\\
 \label{eq:trans2}
   &\uv{y} \times \pq{\g{H}^+ - \g{H}^-} = \g{J}_{\mathrm{s}} = \pt{\underline{\g{Y}}_2 - \underline{\g{Y}}_1} \cdot \g{E}_\tau \,,
\end{align}
where $\pm$ superscripts indicate the quantity evaluated immediately above or below the sheet, it is readily seen that \eqref{eq:sibc2}--\eqref{eq:trans2} 
%
%
imply \eqref{eq:sibc}, thus proving the equivalence of the two configurations.

From \eqref{eq:trans2} and \eqref{eq:emod} the equivalent current density can be written as
\begin{align}
   \g{J}_{\mathrm{s}}\pt{x,z}  &= \g{j}_{\mathrm{s}}\pt{x} \esp{-\jrm \kz z} \,,
\end{align}
where 
\begin{align}
\label{eq:bcmod}
   \g{j}_{\mathrm{s}}\pt{x}  &= \pt{\underline{\g{Y}}_2 - \underline{\g{Y}}_1} \cdot \g{e}_\tau \pt{x,y=0} \,.
\end{align}
In what follows, the spectral (wavenumber) domain is accessed via the Fourier-transform pair
\begin{equation}
\label{eq:ft1}
    \tilde{\g{f}}\pt{\kx} = \intii  \g{f}\pt{x} \esp{+\jrm \kx x} \dd x,
\end{equation}
\begin{equation}
    \label{eq:ft2}
    \g{f}\pt{x} = \fr{1}{2\pi}\intii  \tilde{\g{f}}\pt{\kx} \esp{-\jrm \kx x} \dd\kx \,,
\end{equation}
with the tilde $\tilde{\cdot}$ denoting spectral-domain quantities.
We assume that the sought modal field is transversely evanescent as $x\rightarrow\pm\infty$ so that the equivalent current density in (\ref{eq:bcmod}) is Fourier-transformable. Accordingly, the tangential electric field generated by the current sheet can be expressed as
\begin{align}
\label{eq:e}
    \g{e}_\tau\pt{x,y}&=\nonumber\\
    &\fr{1}{2\pi}\intii \underline{\tilde{\g{G}}}_{Y_1}^{\mathrm{ee}} \pt{\kx,\kz,y,y'=0}\cdot \tilde{\g{j}}_{\mathrm{s}}\pt{\kx} \esp{-\jrm\kx x} \dd\kx  \,,
\end{align}
where $\underline{\tilde{\g{G}}}_{Y_1}^{\mathrm{ee}}\pt{\kx,\kz,y,y'}$ is the spectral dyadic Green's function for the tangential electric field produced by an electric surface current placed at $y'$ above the impedance boundary \eqref{eq:sibc2}; its expression is derived in Appendix \ref{sec:appGreen}.

By inserting \eqref{eq:e} into \eqref{eq:bcmod}, and using \eqref{eq:ft2}, the following integral equation is finally obtained:
\begin{equation}
\label{eq:ie}
  \intii \underline{\tilde{\g{K}}}\pt{\kx,\kz}\cdot \tilde{\g{j}}_{\mathrm{s}}\pt{\kx} \esp{-\jrm \kx x} \dd\kx =0 \,,\quad x>0 \,,
\end{equation}
where 
\begin{equation}
  \underline{\tilde{\g{K}}}\pt{\kx,\kz} = \underline{\g{I}} - \pt{\underline{\g{Y}}_2 - \underline{\g{Y}}_1} \cdot  \underline{\tilde{\g{G}}}_{Y_1}^{\mathrm{ee}} \pt{\kx,\kz,y=0,y'=0}
\end{equation}
and $\underline{\g{I}}$ is the $2\times2$  identity dyadic.

\section{Method-of-Moments Discretization}
\label{Sec:MoM}

\subsection{Entire-Domain Basis Functions}

The Cartesian components of the unknown electric current density $\g{j}_{\mathrm{s}}\pt{x}$ are represented as a linear combinations of basis functions $\Lambda_{\pt{x,z}n}\pt{x}$, i.e.,
\begin{align}
\label{eq:reprx}
    j_{\mathrm{s}x}\pt{x} = \sum_{n=1}^{N_x} I_{xn} \Lambda_{xn}\pt{x} \,, \\
    \label{eq:reprz}
     j_{\mathrm{s}z}\pt{x} = \sum_{n=1}^{N_z} I_{zn} \Lambda_{zn}\pt{x} \,.
\end{align}
In our formulation, we adopt entire-domain basis functions of the type
\begin{equation}\label{eq:bf}
    \Lambda_{\pt{x,z}n}\pt{x} = L_{n-1}\pt{2a\hat{x}} \esp{-a\hat{x}} \,,
\end{equation}
where $L_n(\cdot)$ is the Laguerre polynomial of order $n$, $a$ is a generally complex coefficient with $\re{a}>0$, and the hat $\hat{\cdot}$  indicates normalization with respect to the free-space wavenumber $k_0$, i.e., $\hat{x}=k_0 x$.  Their Fourier transforms are \cite[7.414.6]{BOOK_Gradshtein}
\begin{equation}
\label{eq:bfspect}
    \tilde{\Lambda}_{\pt{x,z}n}\pt{\kx} = \fr{\jrm}{k_0} \fr{\pt{\hat{k}_x-\jrm a}^{n-1}}{\pt{\hat{k}_x+\jrm a}^{n}}
\end{equation}
where $\hat{k}_x=\kx/k_0$. For $a=1/2$, \eqref{eq:bf} constitute a complete and orthogonal set in the space $L^2\left[0,+\ii\right)$ of square-summable functions defined on the positive real axis \cite{shen2001_laguerre}. Nevertheless, by selecting different values for $a$, it may be possible to reduce the number of basis functions required for a given level of accuracy, as will be illustrated in Sec. \ref{sec:results}.

\subsection{Galerkin Testing and MoM Matrix}

By inserting \eqref{eq:reprx}, \eqref{eq:reprz} into \eqref{eq:ie} and performing a Galerkin testing, the integral equation \eqref{eq:ie} is discretized into a homogeneous linear algebraic system that can be written in block form as
\begin{equation}
\label{eq:disc}
\left(\begin{array}{@{}c|c@{}}
  Z^{\mathrm{MoM}}_{xx}
  & Z^{\mathrm{MoM}}_{xz} \\
\hline
  Z^{\mathrm{MoM}}_{zx} &
  Z^{\mathrm{MoM}}_{zz}
\end{array}\right)
\cdot
\begin{pmatrix}
I_x \\
\hline
I_z
\end{pmatrix}
=
\begin{pmatrix}
0 \\
\hline
0
\end{pmatrix} \,,
\end{equation}
where the elements of the MoM matrix are
\begin{equation}
\label{eq:Zmom}
   Z^{\mathrm{MoM}}_{pq,mn} \pt{\kz} = \intii \tilde{\Lambda}_{pm}\pt{-\kx} \tilde{K}_{pq}\pt{\kx,\kz} \tilde{\Lambda}_{qn}\pt{\kx} \dd \kx,
\end{equation}
with $p,q=x,z$ and $m,n=1,2,\dotsc N_{x,z}$. The MoM-matrix elements are thus seen to be functions of the unknown modal wavenumber $\kz$, defined by integrals in the complex plane of the transverse spectral variable $\kx$. 

Note that, by virtue of \eqref{eq:bfspect}, the integrand in \eqref{eq:Zmom} depends only on the difference $\pt{m-n}$, i.e., each of the MoM-matrix block is of Toepliz type. This provides considerable computational advantages, since  the evaluation time of the MoM matrix scales as the number of basis functions rather than its square.

Furthermore, from \eqref{eq:Zmom} and the expression of the spectral dyadic Green's function (see \eqref{eqn:sgfiso} in Appendix \ref{sec:appGreen}), it is also readily established that, when the surface admittances are isotropic (i.e., $\underline{\g{Y}}_{1,2}=Y_{1,2}\underline{\g{I}}$), the diagonal blocks of the MoM matrix are Hermitian, i.e.,
\begin{equation}
\label{eq:symm}
\begin{split}
    \pq{Z^{\mathrm{MoM}}_{xx}}^{\mathrm{T}}=\pq{Z^{\mathrm{MoM}}_{xx}}^* \,, \\
    \pq{Z^{\mathrm{MoM}}_{zz}}^{\mathrm{T}}=\pq{Z^{\mathrm{MoM}}_{zz}}^*
    \end{split}
\end{equation}
and that the MoM matrix is block-anti-symmetric, i.e.,
\begin{equation}
\label{eq:anti}
\pq{Z^{\mathrm{MoM}}_{zx}}=-\pq{Z^{\mathrm{MoM}}_{xz}}^{\mathrm{T}} \,.
\end{equation}
Therefore, the calculation of the MoM matrix requires the numerical evaluation of exactly $N_x+N_z+\max\left\{ N_x,N_z\right\}$  integrals in the complex $\kx$-plane (i.e., the transverse spectral plane).

\subsection{Singularities in the Transverse Spectral Plane}

The evaluation of the MoM-matrix elements requires a careful analysis of the integrands in \eqref{eq:Zmom}, in particular of their singularities in the complex $k_x$-plane.

\begin{figure}[!t]
	\centering
	\includegraphics[width=\columnwidth]{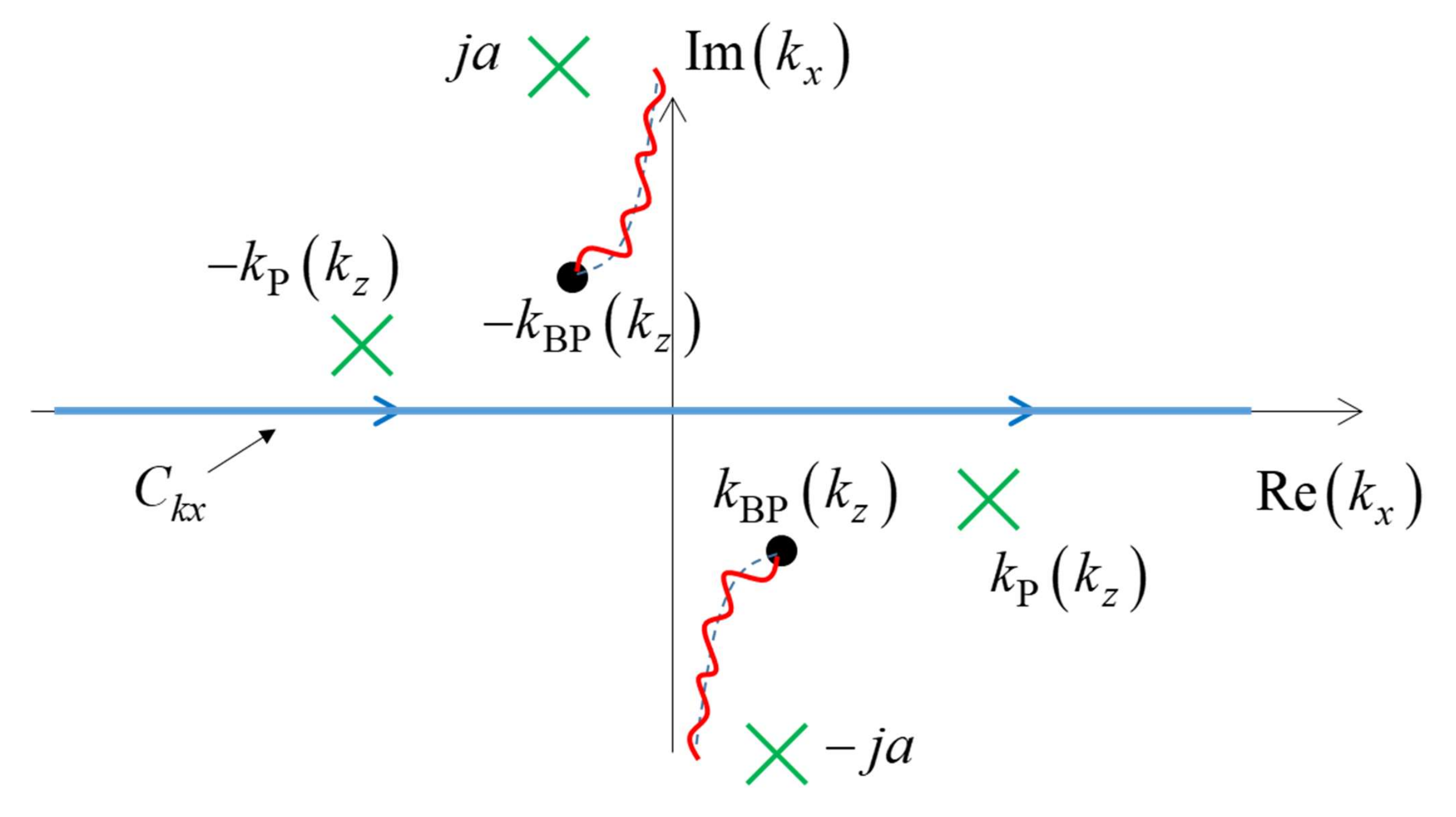}
	\caption{Integration path in the complex $k_x$-plane used to evaluate the elements of the  MoM-matrix for the modal analysis of proper LW modes and the relevant integrand singularities\label{fig:int}}
\end{figure}

First of all, such integrands have a pair of branch-point singularities at $\kx=\pm k_{\mathrm{BP}}\pt{k_z}$, arising from the square-root function that defines the vertical wavenumber $\ky$, with
\begin{equation}
    k_{\mathrm{BP}}\pt{k_z} =  \sqrt{k_0^2-\kz^2} \,.
\end{equation}
The standard hyperbolic Sommerfeld branch cuts defined by $\im{\ky}=0$ then divide the relevant two-sheeted Riemann surface into a proper sheet where $\im{\ky}<0$ and an improper sheet where $\im{\ky}>0$. The integration path to be used for the analysis of transversely evanescent (i.e., proper) modes is thus the real axis of the proper Riemann sheet.

The integrand in \eqref{eq:Zmom} also has pole singularities at $\kx=\pm k_{\mathrm{P}}\pt{k_z}$ which arise from the zeros of the denominator of the Green's function elements (see \eqref{eq:det} in Appendix \ref{sec:appGreen}), corresponding to the SWs supported by a uniform impedance boundary with surface admittance dyadic $\underline{\g{Y}}_1$. Such waves are generally hybrid, since an anisotropic boundary generally couples the transversely magnetic (TM) and electric (TE) polarizations.
%
%
%
%
%
%

Finally, the integrands in \eqref{eq:Zmom} also have a pair of poles at $k_x=\pm \jrm a$, arising from the basis and test functions.

Figure \ref{fig:int} shows the complex $k_x$-plane with the integration path (\textit{blue solid line}) and the integrand singularities: branch points (\textit{black dots}), Sommerfeld branch cuts (\textit{red wiggly lines}), and poles (\textit{green crosses}).

\subsection{Modal Analysis}

The homogeneous system \eqref{eq:disc} defines a non-linear eigenvalue problem, whose solutions can be found by enforcing the condition that the MoM matrix be singular, i.e.,
\begin{equation}
\label{eq:disp}
    \mathrm{det}\pq{Z^{\mathrm{MoM}}\pt{\kz}} = 0 \,,
\end{equation}
which constitutes the dispersion equation for the LW modes.

In our approach, the zeros of \eqref{eq:disp} in the complex $k_z$-plane are computed numerically by using an efficient root-finding procedure based on Padé approximants \cite{galdi2000simple} in which, for each value of $k_z$, the numerical evaluation of the integrals in \eqref{eq:Zmom} is carried out by using a robust and efficient double-exponential quadrature formula \cite{takahasi1974double,michalski2016efficient}.

\section{Numerical Results}\label{sec:results}

In this Section, the proposed spectral MoM formulation is applied to the modal analysis of specific two-part impedance planes characterized by isotropic (Subsec. IV-A) and anisotropic (Subsec. IV-B) boundary conditions. The MoM results are validated against results obtained with Sommerfeld--Maliuzhinets analytical method in the isotropic case, and with finite-element (COMSOL Multiphysics \cite{COMSOL:2015}) simulations in the anisotropic case. 

\subsection{Isotropic Structures}

\subsubsection{Bound Modes}

\begin{figure}[!t]
	\centering
	\includegraphics[width=\columnwidth]{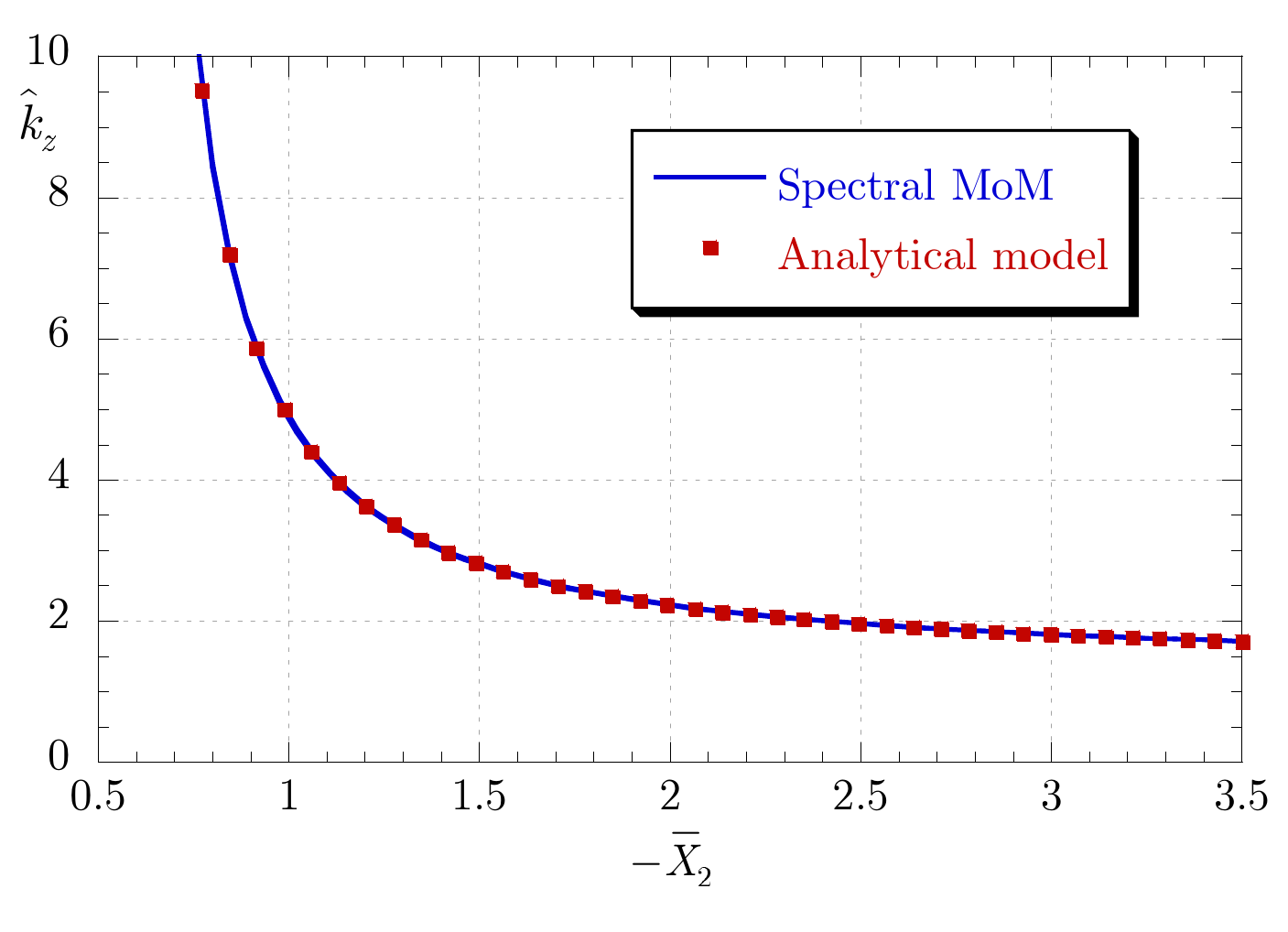}
	\caption{Normalized wavenumber $\hat{k}_z=k_z/k_0$ of the bound LW supported by a complementary reactive two-part impedance plane with ${\bar{Z}_1=\jrm \bar{X}_1=\jrm /\sqrt{3}}$ and $\bar{Z}_2=\jrm \bar{X}_2$ as a function of $-\bar{X}_2$. Comparison between results obtained with the spectral MoM and the Sommerfeld--Maliuzhinets analytical model. \label{fig:disp}}
\end{figure}

Let us consider first a complementary reactive two-part impedance plane, characterized by scalar normalized surface impedances ${\bar{Z}_1=Z_1/\eta_0=\jrm \bar{X}_1}$ and $\bar{Z}_2=Z_1/\eta_0=\jrm \bar{X}_2$ (with $\eta_0$ denoting the free-space characteristic impedance). In Fig. \ref{fig:disp} the normalized wavenumber $\hat{k}_z=k_z/k_0$ of the bound LW supported by such a structure is reported for $\bar{X}_1=1/\sqrt{3}$ as a function of $-\bar{X}_2$ (when $\bar{X}_2=-\sqrt{3}$ the case in \cite[Fig. 3(a)]{kong2019analytic} is recovered); the results obtained with  the spectral MoM and the Sommerfeld--Maliuzhinets analytical model are perfectly superimposed. Note that, as expected \cite{Horsley:2014od}, the LW wavenumber tends to infinity when $-\bar{X}_2$ approaches the value $\bar{X}_1=1/\sqrt{3}\simeq 0.58$.

\begin{figure}[!t]
\centering
    \subfigure[]{
 \includegraphics[width=\columnwidth]{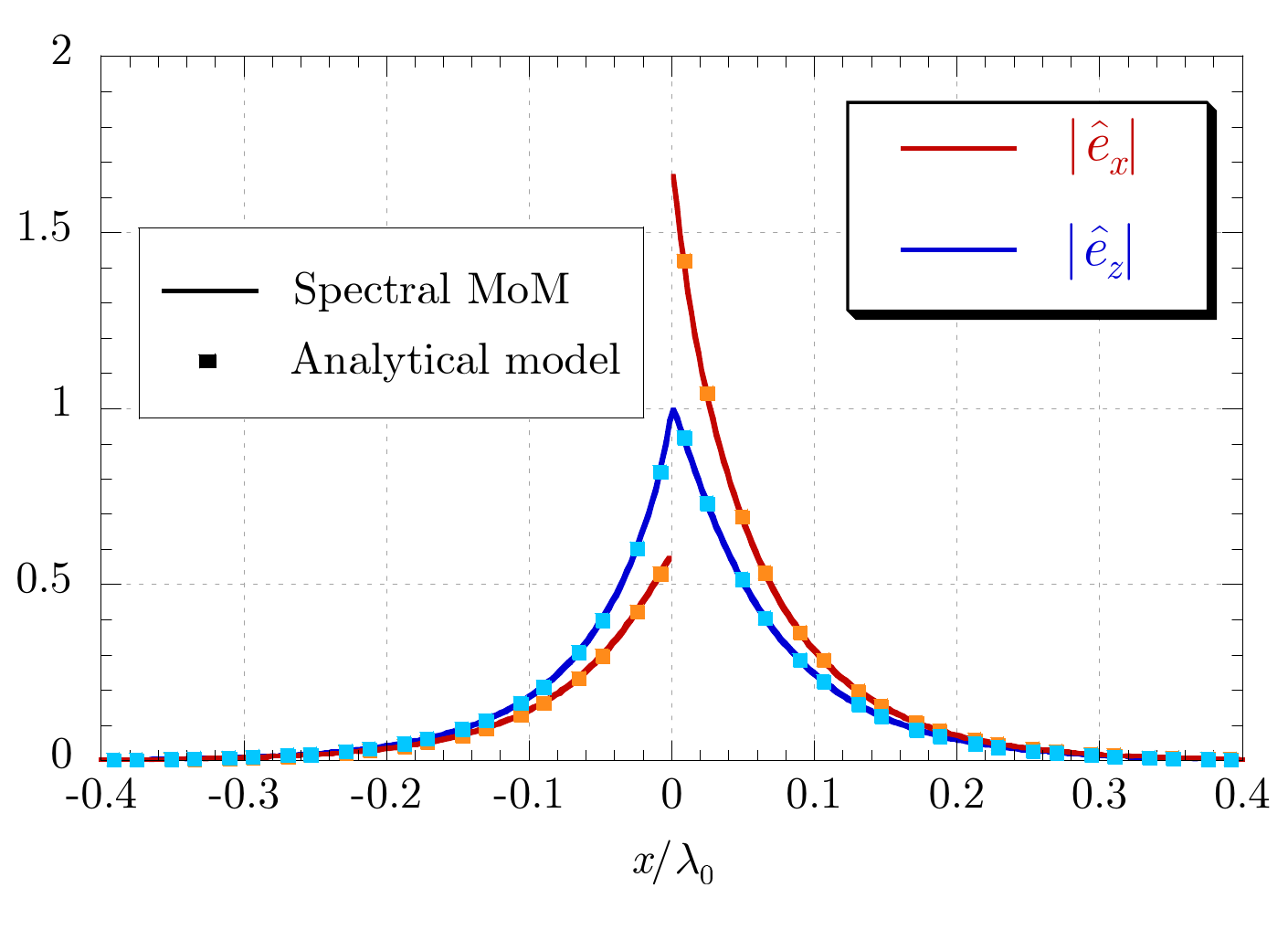}}
 \subfigure[]{
\includegraphics[width=\columnwidth]{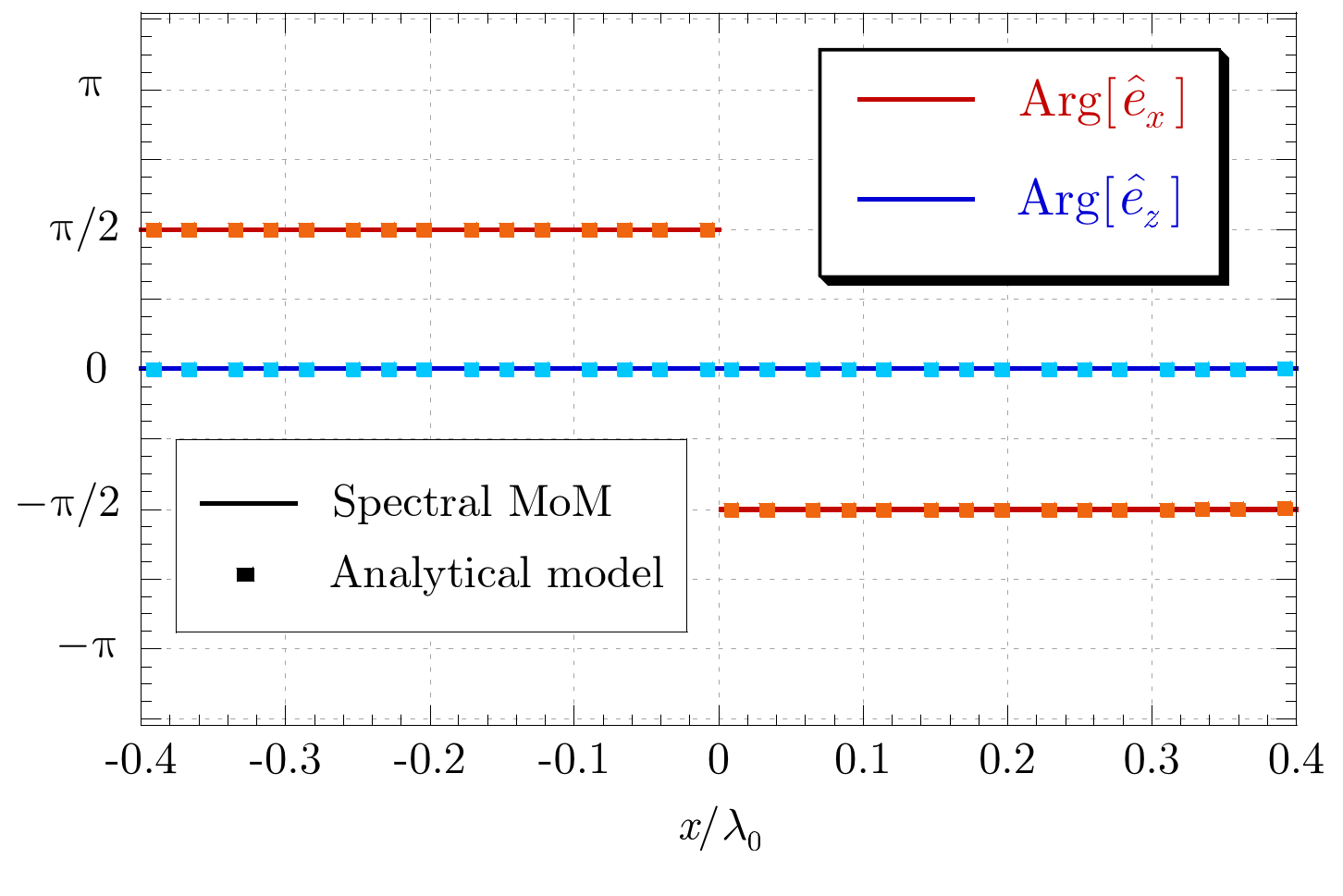}}
	\caption{In-plane components $e_x(x)$ and $e_z(x)$ of the modal field, normalized to $e_z(0)$, for the bound LW supported by a complementary reactive two-part impedance plane with ${\bar{Z}_1=\jrm \bar{X}_1=\jrm /\sqrt{3}}$ and $\bar{Z}_2=\jrm \bar{X}_2=-j\sqrt{3}$, as a function of the normalized transverse coordinate $x/\lambda_0$. Comparison between results obtained with the spectral MoM and the Sommerfeld--Maliuzhinets analytical model. (a) Magnitude; (b) phase.  \label{fig:field}}
\end{figure}

In Fig. \ref{fig:field} the normalized in-plane components ${\hat{e}_x(x)=e_x(x)/e_z(0)}$ and ${\hat{e}_z(x)=e_z(x)/e_z(0)}$ of the relevant modal electric field are reported as a function of the normalized transverse coordinate $x/\lambda_0$ (with $\lambda_0=2\pi/k_0$ denoting the free-space wavelength). Also in this case, the MoM and analytical results are in perfect agreement.

Note that, as expected from the edge conditions derived in \cite{braver2014edge} for a two-part impedance plane, both the transverse and longitudinal components of the field remain finite in the vicinity of $x=0$ in the plane $y=0$. On the other hand, for $y>0$, the radial component $e_\rho$ of the field would exhibit a logarithmic singularity in the vicinity of $\rho=0$, which has been verified numerically both in the spectral MoM and the analytical results (not shown here for brevity). 

It is interesting to observe that the longitudinal component of the modal field $e_z$, parallel to the line discontinuity at $x=0$ between the two impedance half-planes, is \textit{continuous} across such a discontinuity, whereas the transverse component $e_x$ is \textit{discontinuous} both in magnitude and phase; the same behavior (not reported here for brevity) would be observed for $h_z$ and $h_x$, respectively. In fact, the continuity of the field components parallel to the line discontinuity between the two impedance half-planes can be predicted in general via the same standard argument used to establish the continuity of the tangential field components across a surface discontinuity between two bulk media, i.e., applying the Faraday-Neumann-Lenz law to an infinitesimal rectangular contour lying in the plane $y=0$ across the line discontinuity. On the other hand, from the boundary conditions \eqref{eq:sibc} one readily finds that, on the plane $y=0$, ${e_x(x=0^-)=Z_1h_z(x=0^-)}$ and ${e_x(x=0^+)=Z_2h_z(x=0^+)}$; therefore, from ${h_z(x=0^-)=h_z(x=0^+)}$ one concludes that ${e_x(x=0^+)/e_x(x=0^-)=Z_2/Z_1}$,  i.e., the transverse electric field is discontinuous across the line discontinuity, with the ratio between the two limiting values on its two sides being equal to the ratio between the corresponding surface impedances. This can indeed be verified in Fig. \ref{fig:field}, where $Z_2/Z_1=-3$.

Finally, it can also be noted that the two shown in-plane components of the modal field are in phase quadrature and that their phase is constant (or piecewise constant) along $x$.

\begin{figure}[!t]
	\centering
	\includegraphics[width=\columnwidth]{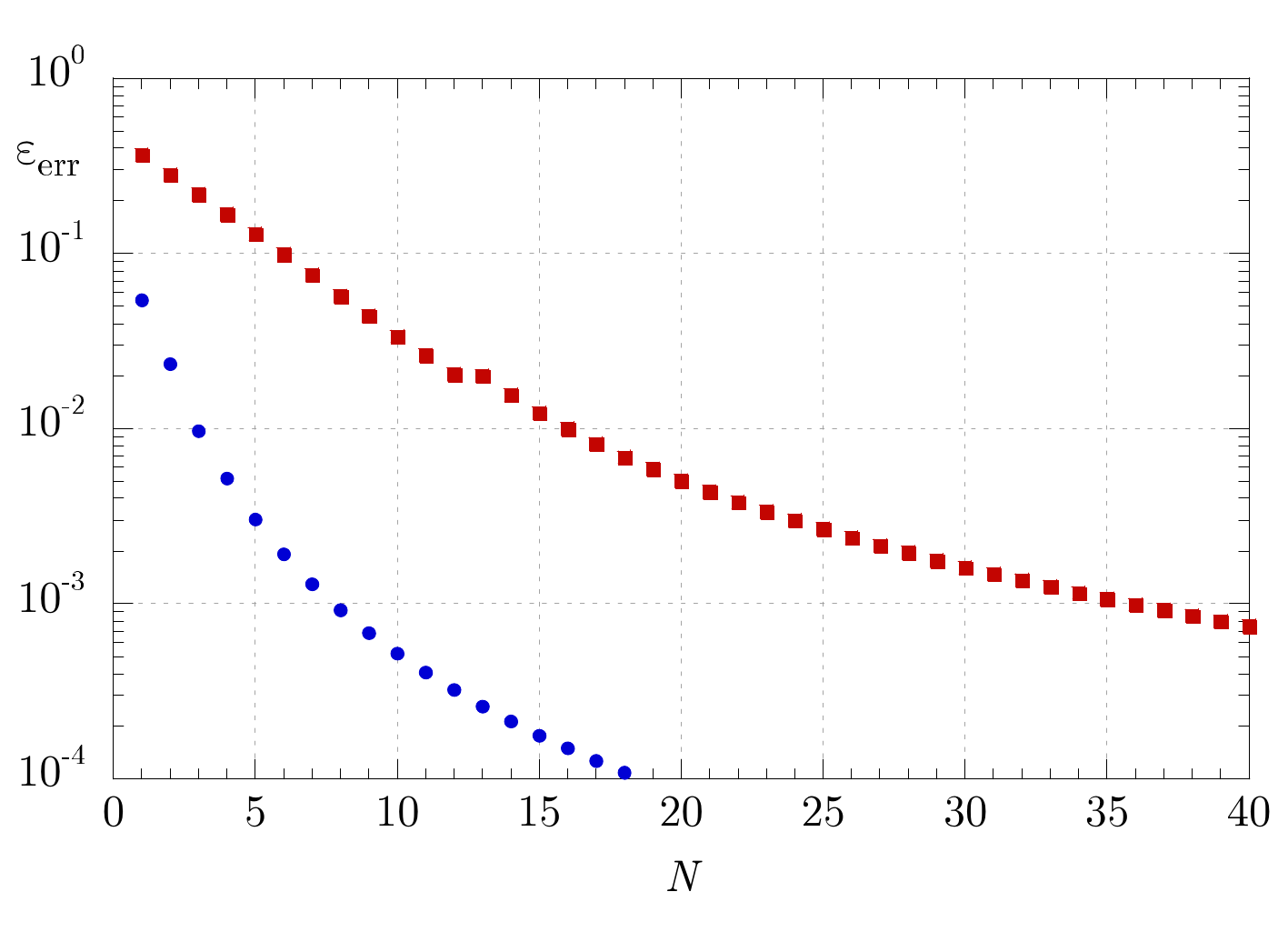}
	\caption{Relative error in the computation of the normalized wavenumber $\hat{k}_z=k_z/k_0$ of the bound LW supported by a structure as in Fig. \ref{fig:field}, as a function of the number $N=N_x=N_z$ of basis functions used to represent each component of the equivalent current density in the MoM formulation, for two choices of the coefficient $a$ in \eqref{eq:bf}: $a=1/2$ (\textit{red squares}) and $a=\sqrt{\hat{k}_z^2-1}$ (\textit{blue circles}). \label{fig:err}}
\end{figure}

Figure \ref{fig:err} shows the relative error in the MoM computation of the modal wavenumber $k_z$ as a function of the number of basis functions $N$ used to represent both the transverse and longitudinal components of the equivalent current density $\g{j}_\mathrm{s}$ (i.e., $N_x=N_z=N$). Two choices are considered for the coefficient $a$ within the exponential function in the basis functions \eqref{eq:bf}, namely the constant $a=1/2$ and the $k_z$-dependent value $a=\sqrt{\hat{k}_z^2-1}$ that would provide the asymptotic radial decay of the field produced by an infinite traveling-wave line source having $z$-dependence $\exp(-jk_zz)$ and parallel to $z$. In both cases the relative error decays monotonically with $N$; however, with the latter choice the decay is considerably faster: for instance, less than $10$ basis functions per component are required to attain an accuracy of $10^{-3}$, and even a single basis function per component provides an accuracy of about $5\%$, sufficient for a rough estimate of the LW wavenumber.

\subsubsection{Leaky Modes}

We now consider a different structure, namely a capacitive two-part impedance plane in which one of the two half-planes is lossy. As shown in \cite{moccia2020line}, this kind of structure may exhibit surface leakage when excited by a three-dimensional source. In fact, the structure supports a \textit{leaky} LW in a surface leakage regime, whose longitudinal wavenumber is therefore complex: $k_z=\beta_z-j\alpha_z$.

\begin{figure}[!t]
	\centering
	\includegraphics[width=\columnwidth]{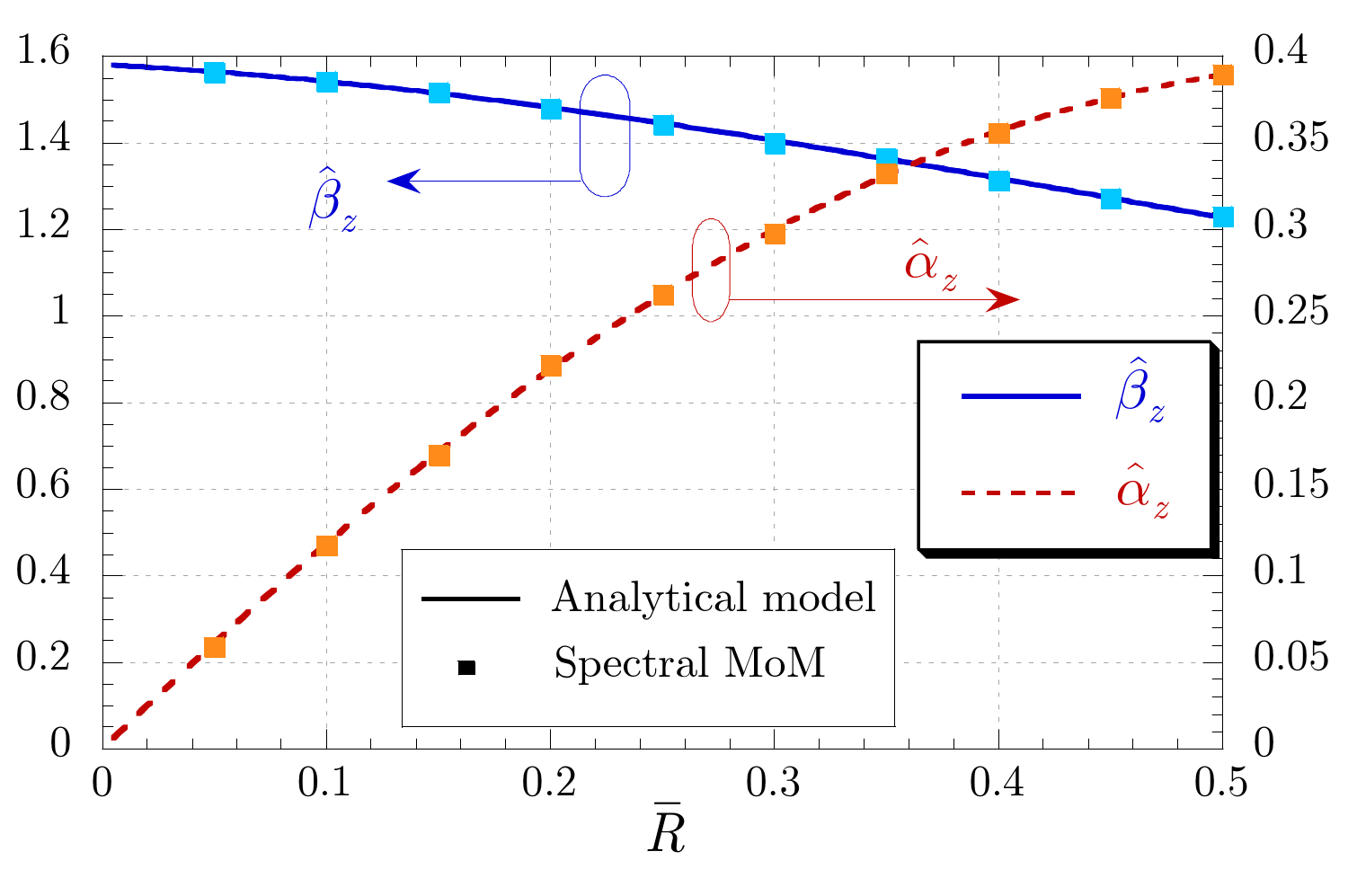}
	\caption{Normalized wavenumber $\hat{k}_z=k_z/k_0$ of the leaky LW supported by a capacitive two-part impedance plane with a lossy half-plane, as a function of the normalized resistance $\bar{R}=R/\eta_0$ of the lossy half-plane. Comparison between results obtained with the spectral MoM and the Sommerfeld--Maliuzhinets analytical model. \textit{Parameters}: ${\bar{Z}_1=-\jrm 0.5}$ and $\bar{Z}_2=\bar{R} -\jrm 0.5$.\label{fig:displeaky}}
\end{figure}

Figure \ref{fig:displeaky} shows the normalized phase constant $\hat{\beta}_z=\beta_z/k_0$ and attenuation constant $\hat{\alpha}_z=\alpha_z/k_0$ of such a leaky mode for a structure with a reactive half-plane ($x<0$) with ${\bar{Z}_1=-\jrm 0.5}$ and a lossy half-plane ($x>0$) with $\bar{Z}_2=\bar{R} -\jrm 0.5$ as a function of the normalized surface resistance $\bar{R}=R/\eta_0$ (when $\bar{R}=0.5$ the case of \cite[Fig. 5(a)]{moccia2020line} is recovered). Within the entire range of $\bar{R}$ values, the normalized phase constant is greater than one but less than the normalized phase constant of the TE SW supported by the reactive half-plane $\beta_\mathrm{SW}^\mathrm{TE}\simeq 2.24$; this reflects the fact that the LW propagates in a \textit{surface leaky} regime, i.e., it leaks power through such a TE SW. Note that, in the limit $\bar{R}\rightarrow 0$, the two-part impedance plane tends to a uniform capacitive plane, and the LW longitudinal wavenumber tends to a real value $\beta_\mathrm{SW}^\mathrm{TE}/\sqrt{2}\simeq 1.6$, corresponding to a leakage angle of $45^\circ$ measured from the longitudinal $z$-axis.

Interestingly, this is a \textit{proper} leaky mode. It has been found with the spectral MoM by selecting the proper determination of the vertical wavenumber $k_y$ in the spectral Green's function $\underline{\tilde{\g{G}}}_{Y_1}^{\mathrm{ee}}$ and by integrating along the real axis of the transverse $k_x$-plane in the evaluation of the MoM-matrix elements. This can be contrasted with the surface leaky modes supported, e.g., by printed-circuit lines, which exhibit an in-plane improper character and require a deformation of the transverse integration path around the poles associated with the SWs through which leakage occurs \cite{mesa2002investigation}.

\begin{figure}[!t]
\centering
    \subfigure[]{
 \includegraphics[width=\columnwidth]{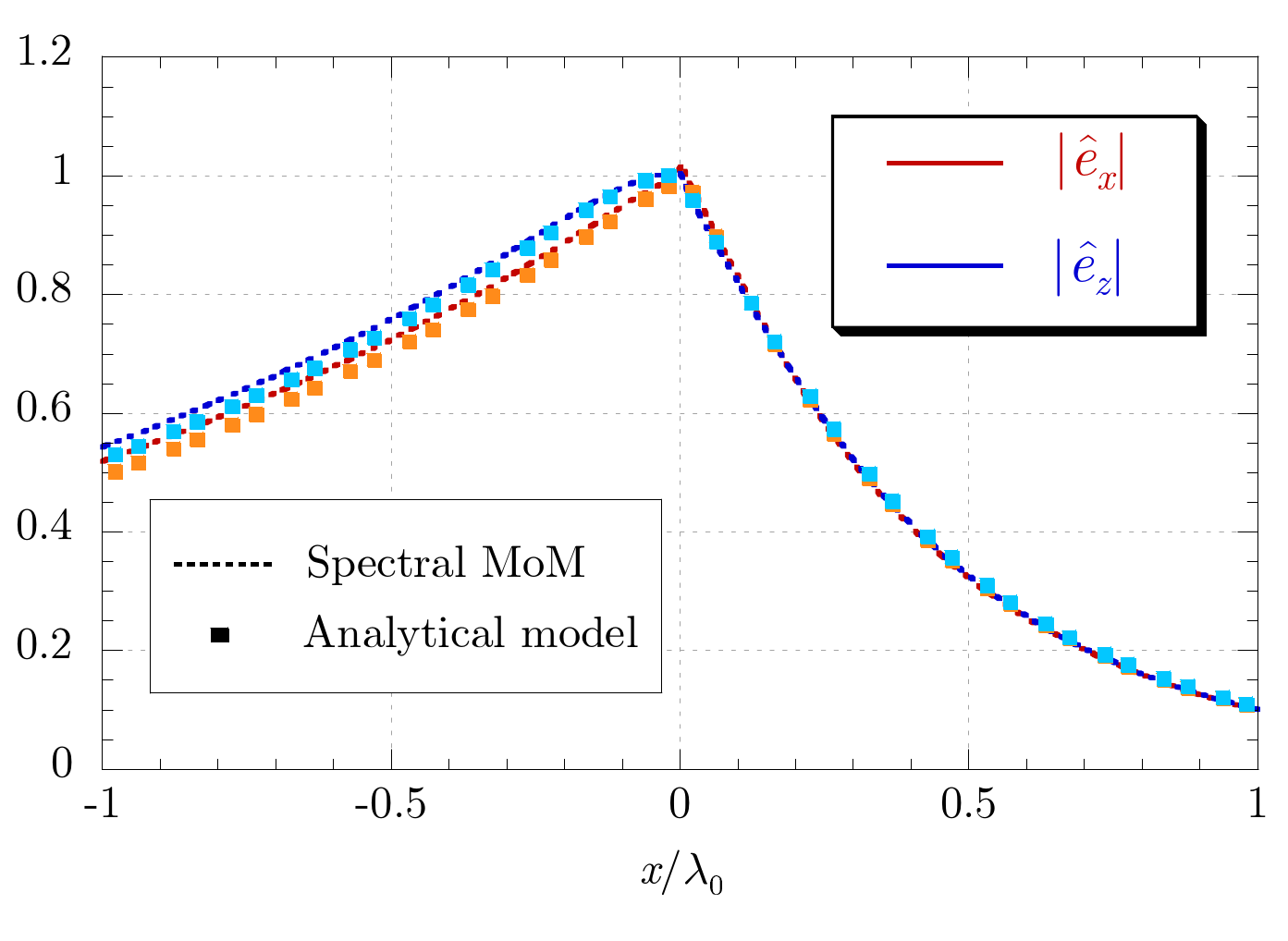}}
 \subfigure[]{
\includegraphics[width=\columnwidth]{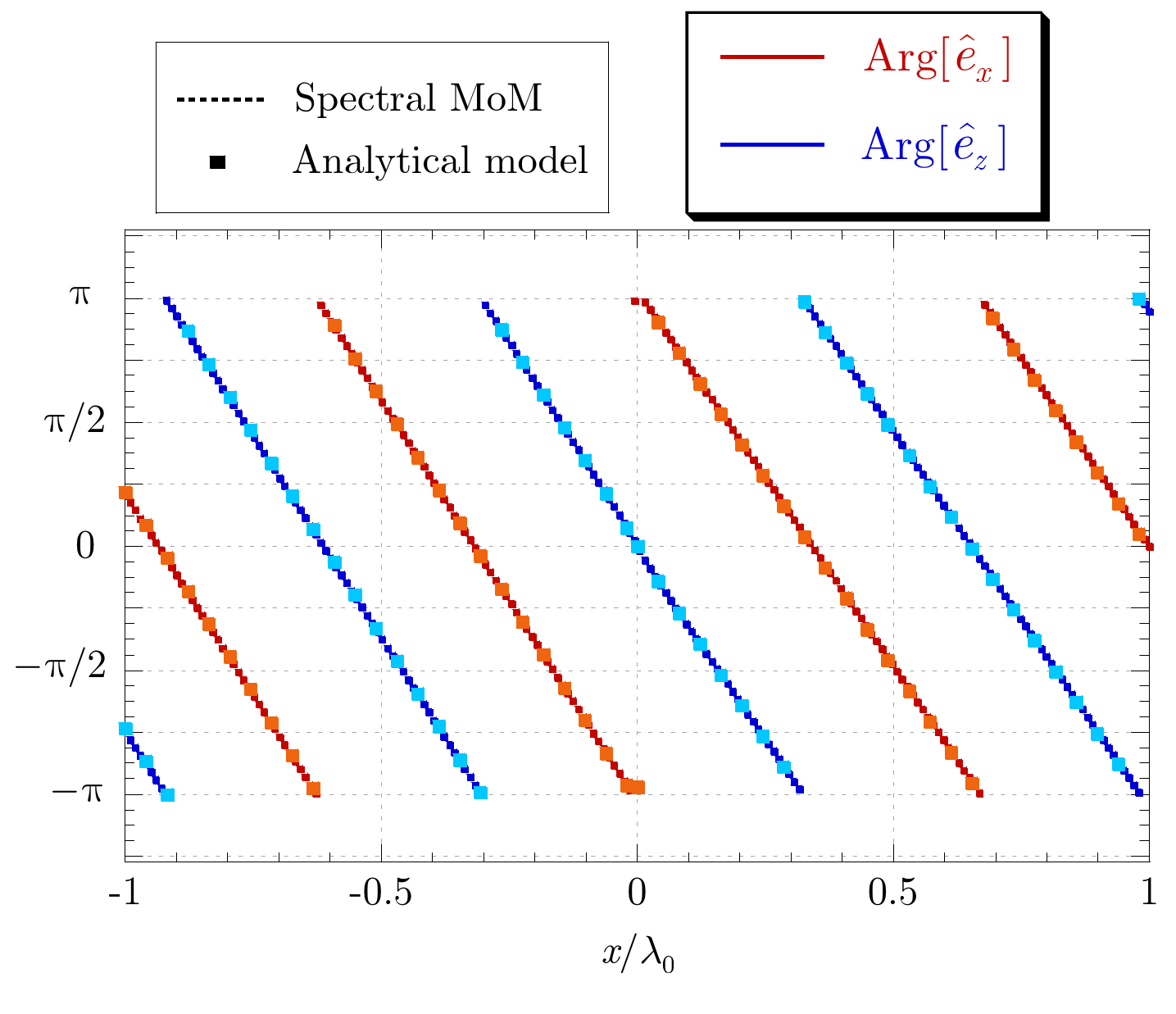}}
	\caption{In-plane components $e_x(x)$ and $e_z(x)$ of the modal field, normalized to $e_z(0)$, for the leaky LW supported by a structure as in Fig. \ref{fig:displeaky} with $\bar{R}=0.1$, as a function of the normalized transverse coordinate $x/\lambda_0$. Comparison between results obtained with the spectral MoM and the Sommerfeld--Maliuzhinets analytical model. (a) Magnitude; (b) phase. \label{fig:fieldleaky}}
\end{figure}

In Fig. \ref{fig:fieldleaky}, the relevant modal electric field is reported for a structure as in Fig. \ref{fig:displeaky} with $\bar{R}=0.1$ (for which $\hat{k}_z \simeq 1.548 - \jrm 0.111$). Once again, the MoM results are in excellent agreement with the analytical results.

In this case, the small impedance contrast between the two half-planes $x<0$ and $x>0$ makes the discontinuity of $e_x$ across $x=0$ negligible. Regarding the field magnitude, it decays exponentially on both sides of the impedance discontinuity, in agreement with the proper nature of the mode, and the transverse decay is seen to be much faster along the lossy half-plane $x>0$ than along the reactive half-plane $x<0$, as expected. In terms of the phase, unlike the bound-mode field in Fig. \ref{fig:field}, both the in-plane components of the modal field exhibit a linear variation of the phase along the transverse coordinate $x$, i.e., a traveling-wave character directed along the positive $x$ axis.  

\begin{figure}[!t]
	\centering
	\includegraphics[width=0.9\columnwidth]{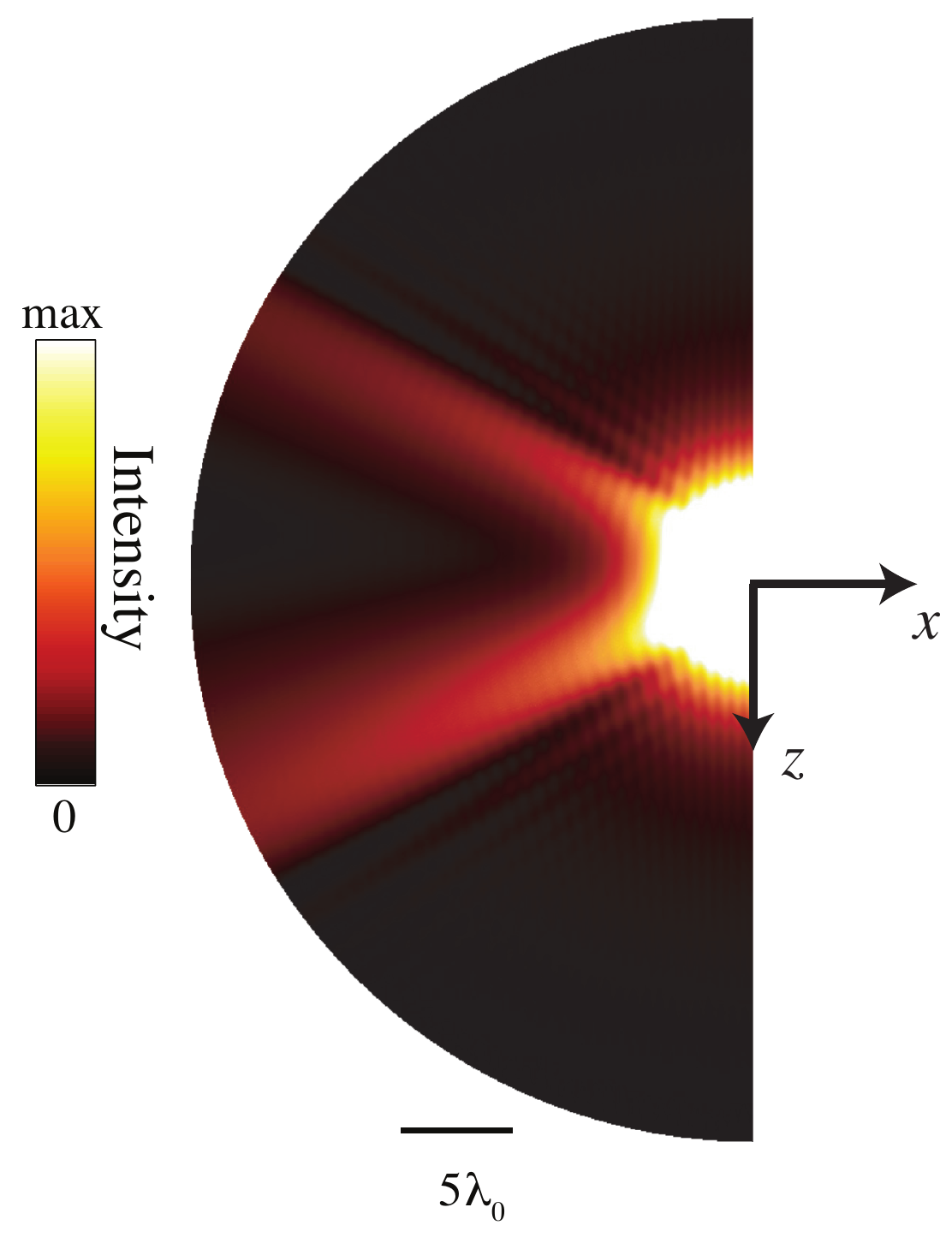}
	\caption{Finite-element computed in-plane field-intensity map ($|{\bf E}|^2$  at $y=0.001\lambda_0$) in false-color scale, for the configuration in Fig. \ref{fig:fieldleaky} (${\bar{Z}_1=-\jrm 0.5}$ and $\bar{Z}_2=0.1 -\jrm 0.5$), illustrating the surface leakage phenomenon. For better visibility, only the lossless region $x<0$ is displayed. The field is excited by a $y$-directed elementary electric dipole placed at $y=0.1\lambda_0$ at the center of the shown semi-circle. The colorscale is suitably saturated so that the space wave radiated by the dipole does not overshadow the surface leakage.\label{fig:mapleaky}}
\end{figure}

The overall modal behavior is reminiscent of the field structure of the Zenneck wave supported by the interface between air and a lossy ground \cite{michalski2015sommerfeld}. In order to illustrate the physical significance of such a proper surface leaky solution, the structure considered in Fig. \ref{fig:fieldleaky} has been excited by a vertical electric dipole placed at the origin of the reference system, i.e., along the line impedance discontinuity. Although  more complex launchers need to be employed in realistic scenarios (see, e.g., \cite{Xu:2019am}), any elementary-dipole excitation is capable in principle of exciting a LW in view of the inherently hybrid character of the relevant field. Specifically, the chosen configuration excites a bi-directional LW propagation (i.e., along the positive and negative $z$ directions), whereas unidirectional propagation can be excited via suitable combinations of dipole sources \cite{dia2017guiding}.
A colormap of the resulting near-field in-plane distribution is shown in Fig. \ref{fig:mapleaky}, from which the typical hallmarks of leaky-wave radiation are observed, with well defined beams of in-plane radiation at angles close to $\cos^{-1}(\beta_z/k_\mathrm{SW})\simeq 46^\circ$ from the $\pm z$ directions.

\subsection{Anisotropic Structures}

\begin{figure}[!t]
	\centering
	\includegraphics[width=\columnwidth]{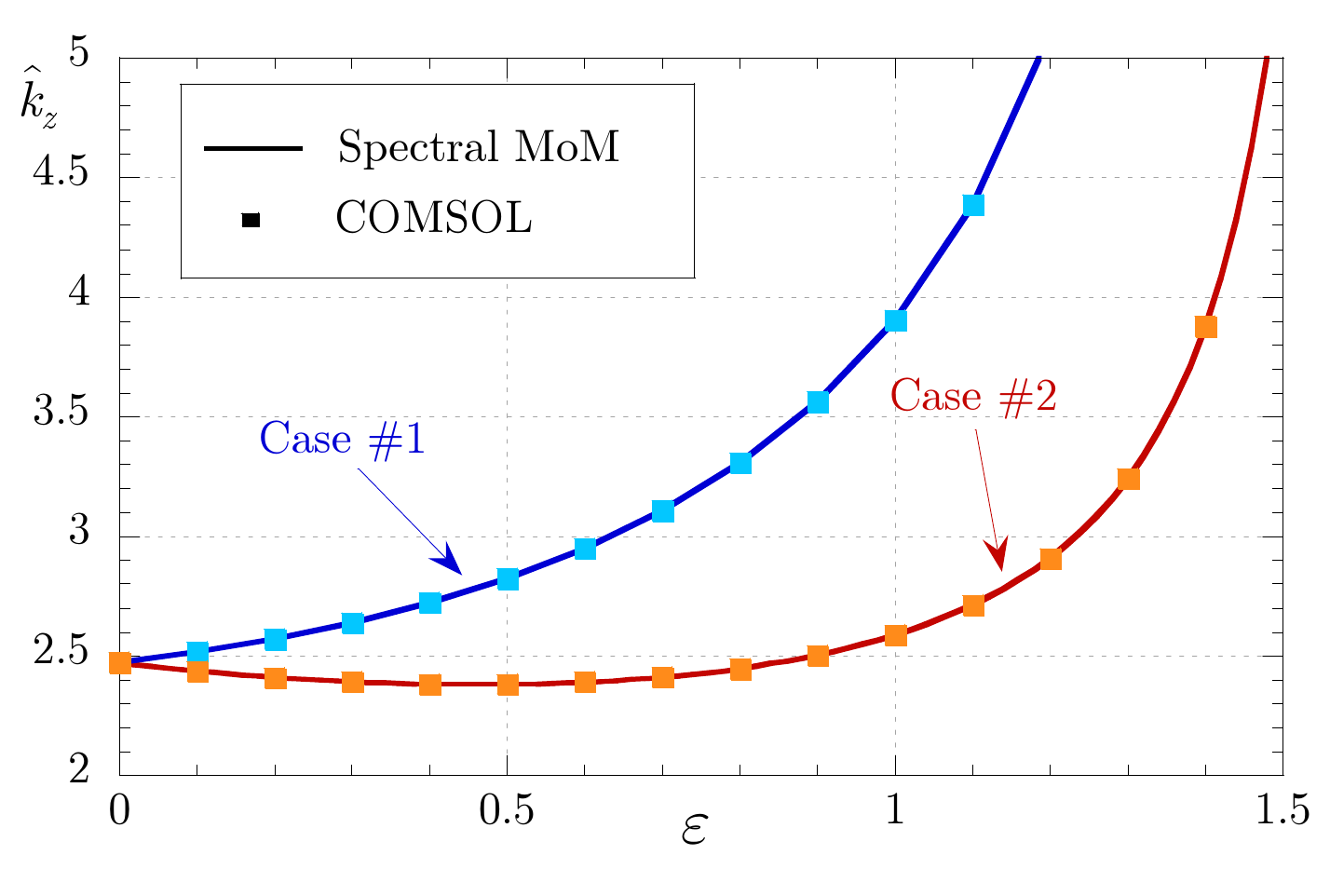}
	\caption{Normalized wavenumber $\hat{k}_z=k_z/k_0$ of the bound LWs supported by complementary reactive two-part impedance planes constituted by an anisotropic half-plane with normalized impedance $\underline{\bar{\g{Z}}}_1$ and an isotropic half-plane with normalized impedance $\bar{Z}_2$, as a function of the anisotropy parameter $\varepsilon$ of $\underline{\bar{\g{Z}}}_1$. Comparison between spectral MoM (\textit{solid lines}) and COMSOL simulations (\textit{squares}). \textit{Parameters}:  ${\underline{\bar{\g{Z}}}_1=-\jrm \left[(\sqrt{3}+\varepsilon)\uz \uz + (\sqrt{3}-\varepsilon)\ux \ux\right]}$ and $\bar{Z}_2=\jrm /\sqrt{3}$  (Case $\#1$: \textit{blue line and cyan squares}); ${\underline{\bar{\g{Z}}}=\jrm \left[1/(\sqrt{3}+\varepsilon)\uz \uz + 1/(\sqrt{3}-\varepsilon)\ux \ux\right]}$ and ${\bar{Z}_2=-\jrm \sqrt{3}}$  (Case $\#2$: \textit{red line and orange squares}). \label{fig:anis1}}
\end{figure}

We now consider LWs supported by structures with anisotropic impedance boundary conditions. In particular, starting from the isotropic structure considered in Fig. \ref{fig:field}, we gradually let the surface impedance of the half-plane $x<0$ become axially anisotropic, considering the following two cases:
\begin{equation}\label{eq:case1}
\begin{split}
&\underline{\bar{\g{Z}}}_1  =-\jrm \left[(\sqrt{3}+\varepsilon)\uz \uz + (\sqrt{3}-\varepsilon)\ux \ux\right]\,, \\
&\bar{Z}_2 = \jrm /\sqrt{3}
\end{split}
\end{equation}
(\textit{Case $\#1$}, in which the anisotropic half-plane $x<0$ is capacitive and the isotropic half-plane $x>0$ is inductive), and
\begin{equation}\label{eq:case2}
\begin{split}
&\underline{\bar{\g{Z}}}_1 = \jrm \left[1/(\sqrt{3}+\varepsilon)\uz \uz + 1/(\sqrt{3}-\varepsilon)\ux \ux\right] \,,\\
&\bar{Z}_2 = -\jrm \sqrt{3}
\end{split}
\end{equation}
(\textit{Case $\#2$}, in which the anisotropic half-plane $x<0$ is inductive and the isotropic half-plane $x>0$ is capacitive). In both cases the real adimensional parameter $\varepsilon$ determines the amount of anisotropy. When $\varepsilon=0$, Case $\#2$ reduces exactly to the isotropic structure of Fig. \ref{fig:field}, whereas Case $\#1$ reduces to the same structure with the surface impedances of the two half planes $x<0$ and $x>0$ interchanged; however, due to the $z$-inversion symmetry of the structure, the modal wavenumber for $\varepsilon=0$ is the same for the two cases, as it can be seen in Fig. \ref{fig:anis1}, where the normalized modal wavenumber $\hat{k}_z$ is reported as a function of $\varepsilon$.

In the axially anisotropic case, an extension of the Sommerfeld--Maliuzhinets method could be developed on the basis of the works available in literature on diffraction by axially anisotropic wedges (see \cite{vallecchi2013oblique} and references therein). However, we have validated here the MoM results against finite-element numerical simulations performed with $\mathrm{{COMSOL}}$ Multiphysics (see Appendix \ref{Sec:DetNum} for details). As shown in Fig. \ref{fig:anis1}, the agreement between the proposed spectral MoM and the finite-element simulations is excellent.

We now consider the value $\varepsilon=1$ for the anisotropy parameter, and rotate the principal anisotropy axes of $\underline{\bar{\g{Z}}}_1$ $\uv{u}$ and $\uv{v}$ of an angle $\xi$:
\begin{equation}
\begin{split}
   &\uv{u} = \cos\xi \uz+\sin\xi \ux \,,\\
  & \uv{v} = -\sin\xi \uz+\cos\xi \ux \,,
\end{split}
\end{equation}
thus considering the surface impedance
\begin{equation}
\underline{\bar{\g{Z}}}_1=-\jrm \left[(\sqrt{3}+1)\uv{u}\uv{u} + (\sqrt{3}-1)\uv{v}\uv{v}\right] \,.
\end{equation}

\begin{figure}[!t]
	\centering
	\includegraphics[width=\columnwidth]{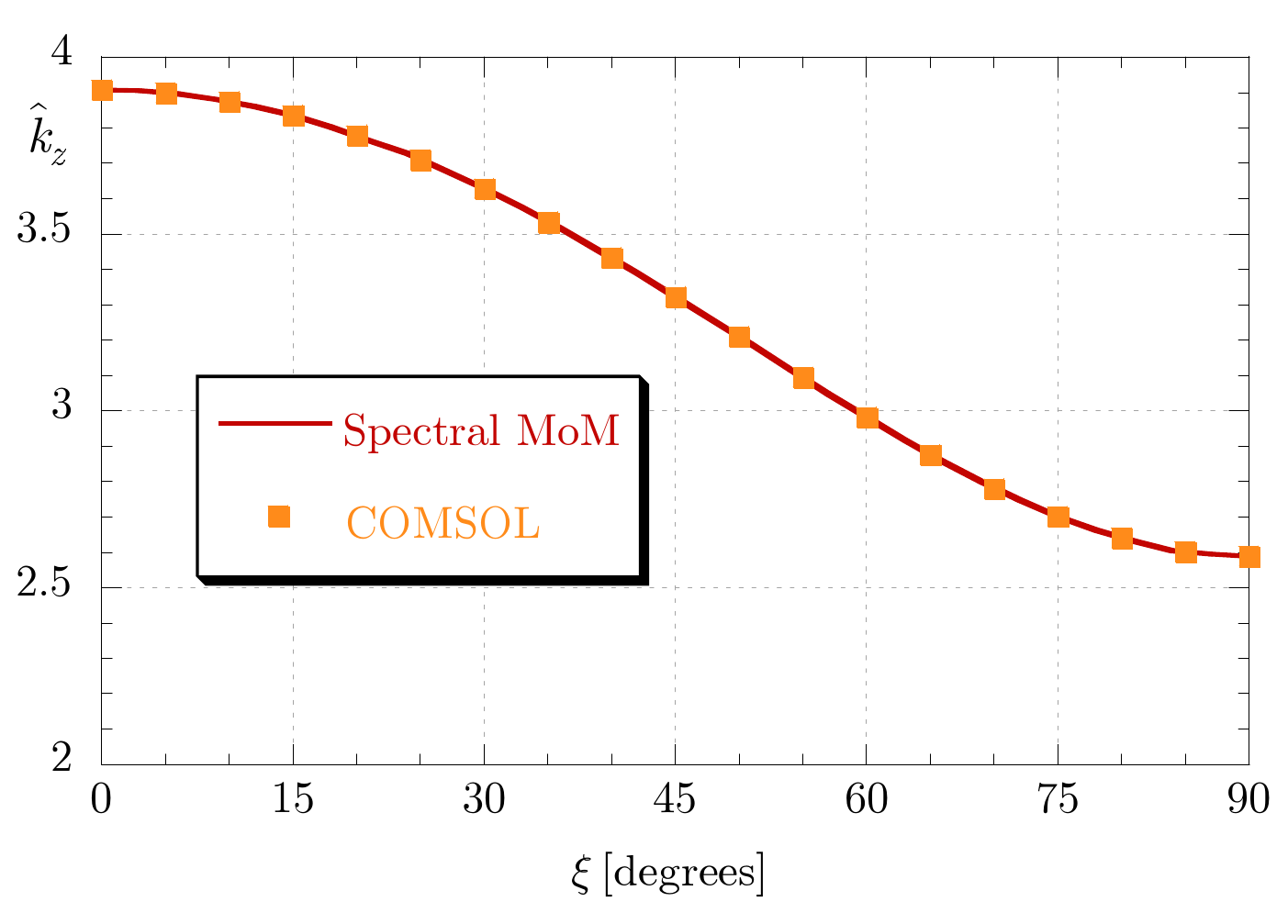}
	\caption{Normalized wavenumber $\hat{k}_z=k_z/k_0$ of the bound LW supported by a complementary reactive two-part impedance plane constituted by an anisotropic half-plane with normalized impedance $\underline{\bar{\g{Z}}}_1$ and an isotropic half-plane with normalized impedance $\bar{Z}_2$, as a function of the skewness angle $\xi$ of the principal anisotropy axes of $\underline{\bar{\g{Z}}}_1$: $\uv{u}=\cos\xi\uz+\sin\xi \ux$ and $\uv{v}=-\sin\xi \uz+\cos\xi \ux$. Comparison between spectral MoM (\textit{solid lines}) and COMSOL simulations (\textit{squares}). \textit{Parameters}: ${\underline{\bar{\g{Z}}}_1=-\jrm \left[(\sqrt{3}+1)\uv{u}\uv{u} + (\sqrt{3}-1)\uv{v}\uv{v}\right]}$ and $\bar{Z}_2=\jrm /\sqrt{3}$. \label{fig:anis2}}
\end{figure}

When $0<\xi<\pi/2$, this surface impedance is \textit{non-axially anisotropic} and no analytical formulation can be adopted for the LW modal analysis (the Sommerfeld--Maliuzhinets approach is available only for particular non-axially anisotropic two-part impedance planes, see \cite{vallecchi2013oblique}). In Fig. \ref{fig:anis2}, the normalized LW modal wavenumber $\hat{k}_z$ is reported as a function of the angle $\xi$, showing a perfect agreement between the MoM and finite-element results. 
Note that the axially anisotropic cases $\#1$ and $\#2$ shown in Fig. \ref{fig:anis1} correspond to specific values of $\xi$. In particular, when $\xi=0$, we recover a specific value of Case $\#1$, and when $\xi=\pi/2$, we recover a specific values of Case $\#2$.

\section{Discussion and Conclusion}
\label{Sec:Conclusions}

The proposed spectral MoM formulation has been shown to be a computationally inexpensive and accurate approach for the modal analysis of line waveguides. It has been validated against independent results considering proper modes supported both by anisotropic structures, for which the Sommerfeld--Maliuzhinets method can be employed, and by general, non-axially anisotropic structures, for which no analytical solution is available.

The proposed method is also versatile, with potential for generalization in various directions. For instance, it can be used to study \textit{improper modes}, by performing the integrations in the transverse spectral plane $k_x$, required to calculate the elements of the MoM matrix, partially or totally on the improper $k_x$-Riemann sheet; this may involve deforming the integration path off the real axis, by suitably adapting the procedure adopted for the modal analysis of printed-circuit lines (see, e.g., \cite{mesa2002investigation} and references therein). 

On the other hand, since the MoM is formulated in the spectral domain, \textit{nonlocal} (i.e., spatially dispersive) impedance boundaries can readily be accommodated by letting the elements of the relevant dyadic surface admittances be dependent on the wavenumbers $(k_x,k_z)$. In this case, however, attention should be paid to the convergence of the resulting spectral integrals, which can be adversely affected by such wavenumber dependence and may even require the adoption of a different testing scheme in the MoM formulation \cite{lovat2015nonlocal}. Within this framework, it appears interesting to compare the MoM predictions with independent numerical simulations of more realistic configurations featuring patterned metallic sheets.

\begin{figure}[!t]
	\centering
	\includegraphics[width=\columnwidth]{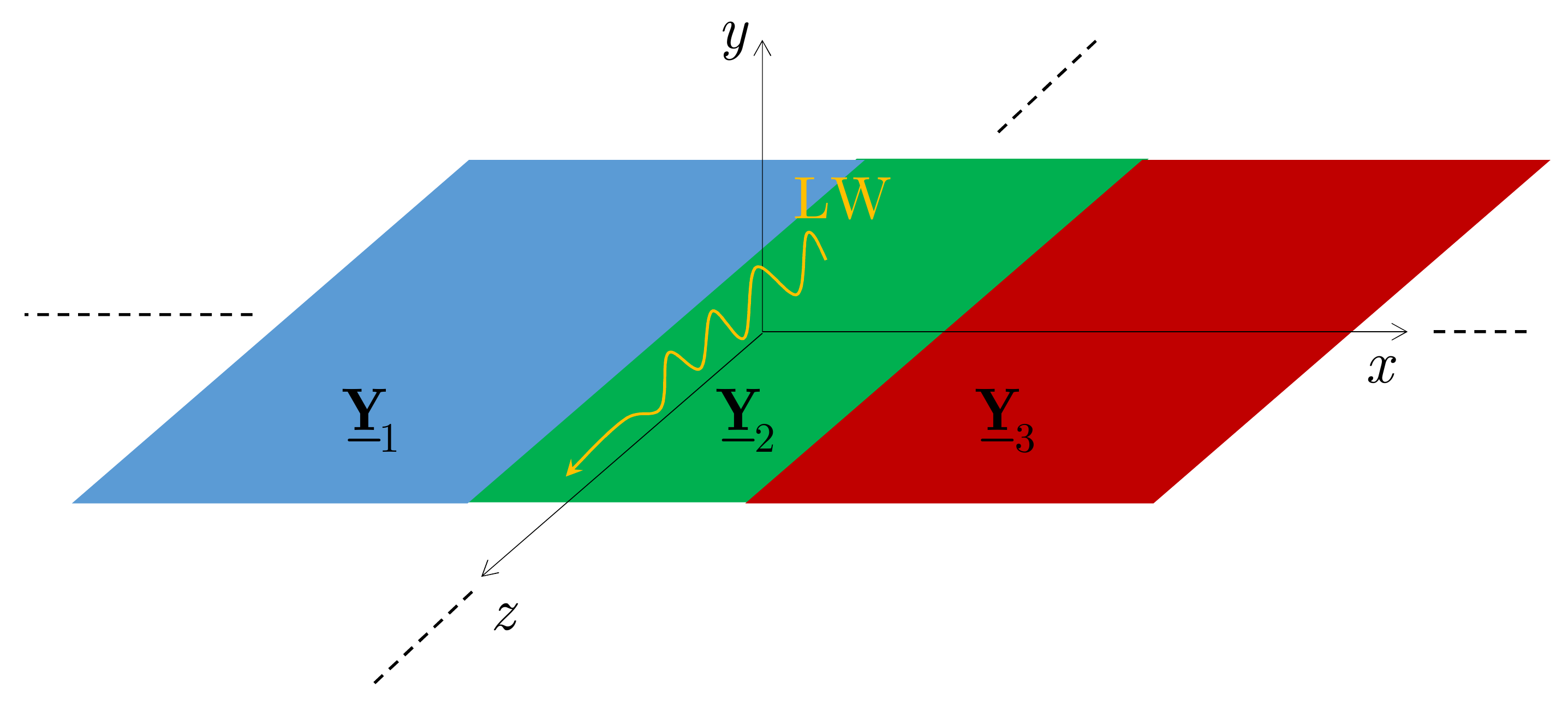}
	\caption{A three-part impedance plane, constituted by two half-planes with dyadic surface admittances $\underline{\g{Y}}_{1,3}$, and a finite-width strip with admittance $\underline{\g{Y}}_{2}$.} \label{fig:struct_3part} 
\end{figure}

In a still different direction, the MoM approach can be extended to deal with impedance planes characterized by \textit{more than one line discontinuity} of their surface impedance. Consider, for instance, the three-part impedance plane shown in Fig. \ref{fig:struct_3part}, constituted by two half-planes with surface admittances $\underline{\g{Y}}_{1,3}$ separated by a finite-width strip with surface admittance. $\underline{\g{Y}}_{2}$. This structure, which is a sort of 2-D version of a dielectric rib waveguide when $\underline{\g{Y}}_{1}=\underline{\g{Y}}_{3}$, has recently been shown to support interesting LW modal regimes \cite{Moccia:2021ep}. It can be studied with our spectral MoM by extending to the entire plane $y=0$ the boundary conditions that holds on one of the two half-planes, as described in Sec. III, and defining two sets of basis functions to represent the equivalent currents placed above the finite-width strip and the other half-plane, respectively.

Work is ongoing to explore these and other possible generalizations.

\appendices

\section{Spectral Green's Function}\label{sec:appGreen}

In this Appendix, we derive the spectral Green's function ${\tilde{\g{G}}^{\mathrm{ee}}_{Y_1}\pt{k_z,k_x;y,y'}}$ used to evaluate the spectral tangential electric field $\tilde{\g{E}}\pt{k_z,k_x,y}$ produced by an electric current sheet $\tilde{\g{J}}_{\mathrm{s}}\pt{k_z,k_x}$ placed at $y'=0$ above the anisotropic planar boundary $y=0$ with spectral dyadic surface admittance
\begin{equation}
    \underline{\g{Y}}_1 = Y_{1zz} \uz \uz + Y_{1zz} \uz \ux + Y_{1xz} \ux \uz + Y_{1xx} \ux \ux \,.
\end{equation}

To this aim, we project fields and currents onto the spectral TM/TE basis $\pt{\uv{u},\uv{v}}$, where 
\begin{equation}
\begin{split}
&\uv{u}=\fr{k_z}{k_\trm}\uz+\fr{k_x}{k_\trm}\ux \,, \\
&\uv{v}=\uy\times\uv{u} \,,
\end{split}
\end{equation}
and $k_\trm=\sqrt{k_z^2+k_x^2}$, and write the boundary condition \eqref{eq:sibc2} as
\begin{equation}\label{eq:spectBC}
\begin{pmatrix}
-\tilde{H}_v \\
\tilde{H}_u
\end{pmatrix}
=
\begin{pmatrix}
Y_{1uu} & Y_{1uv} \\
Y_{1vu} & Y_{1vv}
\end{pmatrix}
\begin{pmatrix}
\tilde{E}_u \\
\tilde{E}_v
\end{pmatrix} \,,
\end{equation}
where $Y_{1pq}=\uv{p}\cdot\underline{\g{Y}}_1\cdot\uv{q}$ with $\pt{p,q}=\pt{u,v}$ or, explicitly,
\begin{equation}
\begin{split}
&Y_{1uu} = \frac{1}{k_t^2}
 \pq{k_z^2 Y_{1zz} + k_z k_x \pt{Y_{1zx} + Y_{1xz}} + k_x^2 Y_{1xx}} \,, \\
&Y_{1uv} = \frac{1}{k_t^2}
 \pq{k_z^2 Y_{1zx} + k_z k_x \pt{Y_{1xx} - Y_{1zz}} - k_x^2 Y_{1xz}} \,,\\
&Y_{1vu} = \frac{1}{k_t^2}
 \pq{k_z^2 Y_{1xz} + k_z k_x \pt{Y_{1xx} - Y_{1zz}} - k_x^2 Y_{1zx}} \,,\\
&Y_{1vv} = \frac{1}{k_t^2}
 \pq{k_z^2 Y_{1xx} - k_z k_x \pt{Y_{1zx} + Y_{1xz}} + k_x^2 Y_{1zz}} \,.
\end{split}
\end{equation}

By introducing the standard equivalent transmission lines associated with TM and TE waves \cite{michalski1997multilayered}, we let 
\begin{equation}
\begin{split}
&V^{\mathrm{TM}}=\tilde{E}_u\pt{k_z,k_x,0} \,, \\
&I^{\mathrm{TM}}=\tilde{H}_v\pt{k_z,k_x,0} \,, \\
&V^{\mathrm{TE}}=\tilde{E}_v\pt{k_z,k_x,0} \,, \\ 
&I^{\mathrm{TE}}=-\tilde{H}_u\pt{k_z,k_x,0} \,, 
\end{split}
\end{equation}
so that  \eqref{eq:spectBC} reads
\begin{equation}\label{eq:viBC2}
\begin{pmatrix}
-I^{\mathrm{TM}} \\
- I^{\mathrm{TE}}
\end{pmatrix}
=
\begin{pmatrix}
Y_{1uu} & Y_{1uv} \\
Y_{1vu} & Y_{1vv}
\end{pmatrix}
\begin{pmatrix}
V^{\mathrm{TM}} \\
V^{\mathrm{TE}}
\end{pmatrix} \,.
\end{equation}

\begin{figure}[!t]
	\centering
	\includegraphics[width=\columnwidth]{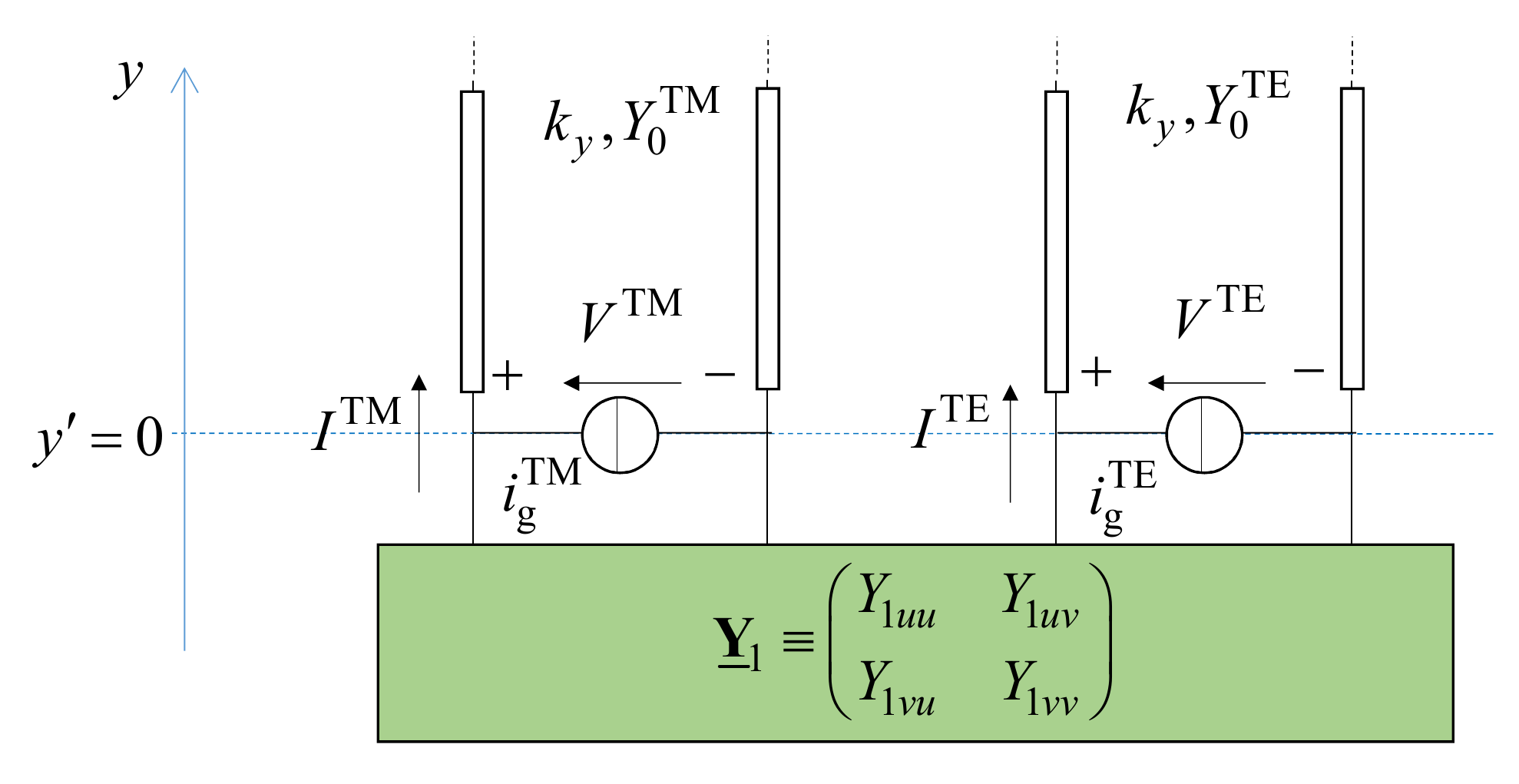}
	\caption{Transverse equivalent network (TEN) for the evaluation of the spectral dyadic Green's function.} \label{fig:TEN}
\end{figure}

The transverse equivalent network (TEN) for the configuration under analysis is thus the one shown in Fig. \ref{fig:TEN}, where the vertical wavenumber in air is $k_y=\sqrt{k_0^2-k_\trm^2}$, the TM/TE characteristic admittances are $Y_0^\mathrm{TM}=k_0/\pt{\eta_0 k_y}$, $Y_0^\mathrm{TE}=k_y/\pt{\eta_0 k_0}$ and the amplitudes of the equivalent current generators are $i_{\mathrm{g}}^{\mathrm{TM}}=-\tilde{J}_u$, $i_{\mathrm{g}}^{\mathrm{TE}}=-\tilde{J}_v$.

The analysis of the TEN is straightforward and leads to
\begin{equation}
\begin{pmatrix}
Y_{1uu}+Y_0^\mathrm{TM} & Y_{1uv} \\
Y_{1vu} & Y_{1vv}+Y_0^\mathrm{TE}
\end{pmatrix}
\begin{pmatrix}
V^\mathrm{TM}\\
V^\mathrm{TE}
\end{pmatrix}
=
\begin{pmatrix}
i_\mathrm{g}^\mathrm{TM}\\
i_\mathrm{g}^\mathrm{TE}
\end{pmatrix},
\end{equation}
from which
\begin{equation}
\begin{split}
\begin{pmatrix}
V^\mathrm{TM}\\
V^\mathrm{TE}
\end{pmatrix}
&=
\begin{pmatrix}
Y_{1uu}+Y_0^\mathrm{TM} & Y_{1uv} \\
Y_{1vu} & Y_{1vv}+Y_0^\mathrm{TE}
\end{pmatrix}^{-1}
\begin{pmatrix}
i_\mathrm{g}^\mathrm{TM}\\
i_\mathrm{g}^\mathrm{TE}
\end{pmatrix}
\\
&=
\begin{pmatrix}
V_{i^\mathrm{TM}}^\mathrm{TM} & V_{i^\mathrm{TE}}^\mathrm{TM} \\
V_{i^\mathrm{TE}}^\mathrm{TM} & V_{i^\mathrm{TE}}^\mathrm{TE}
\end{pmatrix}
\begin{pmatrix}
i_\mathrm{g}^\mathrm{TM}\\
i_\mathrm{g}^\mathrm{TE}
\end{pmatrix} \,.
\end{split}
\end{equation}

The latter expression defines the generalized network Green's functions
\begin{equation}
\begin{split}
&V_{i^\mathrm{TM}}^\mathrm{TM} = \frac{Y_{1vv}+Y_0^\mathrm{TE}}{\mathrm{det}} \,,\\
&V_{i^\mathrm{TE}}^\mathrm{TM} = -\frac{Y_{1uv}}{\mathrm{det}} \,,\\
&V_{i^\mathrm{TM}}^\mathrm{TE} = -\frac{Y_{1vu}}{\mathrm{det}} \,,\\ &V_{i^\mathrm{TE}}^\mathrm{TE} = \frac{Y_{1uu}+Y_0^\mathrm{TM}}{\mathrm{det}} \,,
\end{split}
\end{equation}
where
\begin{equation}\label{eq:det}
\mathrm{det} = \pt{Y_{1uu}+Y_0^\mathrm{TM}}\pt{Y_{1vv}+Y_0^\mathrm{TE}} - Y_{1uv} Y_{1vu} \,.
\end{equation}

We may thus write
\begin{equation}
\begin{split}
\begin{pmatrix}
\tilde{E}_u\\
\tilde{E}_v
\end{pmatrix}
=
\begin{pmatrix}
V^\mathrm{TM}\\
V^\mathrm{TE}
\end{pmatrix}
&=
\begin{pmatrix}
V_{i^\mathrm{TM}}^\mathrm{TM} & V_{i^\mathrm{TE}}^\mathrm{TM} \\
V_{i^\mathrm{TM}}^\mathrm{TE} & V_{i^\mathrm{TE}}^\mathrm{TE}
\end{pmatrix}
\begin{pmatrix}
i_\mathrm{g}^\mathrm{TM}\\
i_\mathrm{g}^\mathrm{TE}
\end{pmatrix} \\
&=
-
\begin{pmatrix}
V_{i^\mathrm{TM}}^\mathrm{TM} & V_{i^\mathrm{TE}}^\mathrm{TM} \\
V_{i^\mathrm{TM}}^\mathrm{TE} & V_{i^\mathrm{TE}}^\mathrm{TE}
\end{pmatrix}
\begin{pmatrix}
\tilde{J}_u\\
\tilde{J}_v
\end{pmatrix}
\end{split}
\end{equation}
and finally, reverting back to Cartesian coordinates, 
\begin{equation}
\begin{pmatrix}
\tilde{E}_z\\
\tilde{E}_x
\end{pmatrix}
=
\begin{pmatrix}
\tilde{G}_{Y_1 zz}^\mathrm{ee} & \tilde{G}_{Y_1 zx}^\mathrm{ee} \\
\tilde{G}_{Y_1 xz}^\mathrm{ee} & \tilde{G}_{Y_1 xx}^\mathrm{ee}
\end{pmatrix}
\begin{pmatrix}
\tilde{J}_z\\
\tilde{J}_x
\end{pmatrix} \,,
\end{equation}
where
\begin{equation}\label{eqn:sgfan}
\begin{split}
   \tilde{G}_{Y_1 zz}^\mathrm{ee} = -\frac{1}{k_{\trm}^{2}}
 \pq{k_z^2 V_{i^\mathrm{TM}}^\mathrm{TM} - k_z k_x \pt{V_{i^\mathrm{TE}}^\mathrm{TM} + V_{i^\mathrm{TM}}^\mathrm{TE}} + k_x^2 V_{i^\mathrm{TE}}^\mathrm{TE}} \,,\\
   \tilde{G}_{Y_1 zx}^\mathrm{ee} = -\frac{1}{k_{\trm}^{2}}
 \pq{k_z^2 V_{i^\mathrm{TE}}^\mathrm{TM} + k_z k_x \pt{V_{i^\mathrm{TM}}^\mathrm{TM} - V_{i^\mathrm{TE}}^\mathrm{TE}} - k_x^2 V_{i^\mathrm{TM}}^\mathrm{TE}} \,,\\
   \tilde{G}_{Y_1 xz}^\mathrm{ee} = -\frac{1}{k_{\trm}^{2}}
 \pq{k_z^2 V_{i^\mathrm{TM}}^\mathrm{TE} + k_z k_x \pt{V_{i^\mathrm{TM}}^\mathrm{TM} - V_{i^\mathrm{TE}}^\mathrm{TE}} - k_x^2 V_{i^\mathrm{TE}}^\mathrm{TM}} \,,\\
   \tilde{G}_{Y_1 xx}^\mathrm{ee} = -\frac{1}{k_{\trm}^{2}}
 \pq{k_x^2 V_{i^\mathrm{TM}}^\mathrm{TM} + k_z k_x \pt{V_{i^\mathrm{TE}}^\mathrm{TM} + V_{i^\mathrm{TM}}^\mathrm{TE}} + k_z^2 V_{i^\mathrm{TE}}^\mathrm{TE}} \,.
 \end{split}
\end{equation}
In the isotropic case, $V_{i^\mathrm{TE}}^\mathrm{TM}=V_{i^\mathrm{TM}}^\mathrm{TE}=0$ and the spectral Green's function reduces to the canonical form (see \cite{michalski1997multilayered})
\begin{equation}\label{eqn:sgfiso}
\begin{split}
   &\tilde{G}_{Y_1 zz}^\mathrm{ee} = -\frac{1}{k_{\trm}^{2}}
 \pt{k_z^2 V_{i^\mathrm{TM}}^\mathrm{TM}  + k_x^2 V_{i^\mathrm{TE}}^\mathrm{TE}} \,, \\
   &\tilde{G}_{Y_1 zx}^\mathrm{ee} = \tilde{G}_{Y_1 xz}^\mathrm{ee}= \frac{k_z k_x}{k_{\trm}^{2}}
   \pt{V_{i^\mathrm{TE}}^\mathrm{TE}-V_{i^\mathrm{TM}}^\mathrm{TM} }  \,,\\
   &\tilde{G}_{Y_1 xx}^\mathrm{ee} = -\frac{1}{k_{\trm}^{2}}
 \pt{k_x^2 V_{i^\mathrm{TM}}^\mathrm{TM}  + k_z^2 V_{i^\mathrm{TE}}^\mathrm{TE}} \,.
 \end{split}
\end{equation}

\section{Details on Finite-Element Simulations}
\label{Sec:DetNum}
The finite-element numerical simulations utilized in Sec. \ref{sec:results} for validation purposes are carried out by means of the commercial software package COMSOL Multiphysics \cite{COMSOL:2015}. In particular, for the eigenmode analysis (Figs. \ref{fig:anis1} and \ref{fig:anis2}), a 2-D computational domain is assumed, comprising a semicircle of radius $2\lambda_0$, with  a perfectly matched layer (PML) termination of thickness $0.5\lambda_0$, and the impedance boundary condition at the plane $y = 0$ enforced via a surface current density ${\bf J}_s={\underline {\bf Y}}\cdot{\bf E}_\tau$ (with ${\underline {\bf Y}}$ particularized for the two halves of the junction). The domain is discretized via an adaptive mesh with element size $\le 0.01\lambda_0$ (refined up to $\le 0.001\lambda_0$ in the vicinity of the impedance surface), which yields $\sim 3$ million degrees of freedom. 
The problem is finally solved by means of the MUMPS direct solver (with default parameters) available in the ``Mode Analysis'' study of the RF Module.

For the study of the dipole-excited configuration in Fig. \ref{fig:mapleaky}, a 3-D computational domain is assumed with out-of-plane height (along $y$) of $0.5\lambda_0$ composed of a semi-cylinder of in-plane radius of $34\lambda_0$  (containing the lossless region $x<0$) and a parallelepiped of in-plane size $34\lambda_0\times0.5\lambda_0$  (containing the lossy region $x>0$), with a  $0.5\lambda_0$-thick PML termination, and the same impedance boundary condition at $y=0$. 
To mimic an infinite impedance surface, fictitious sections of length $9\lambda_0$ are added at the ends, with a linearly tapered profile of the surface impedance that eventually matches the free-space characteristic impedance at the PML. An adaptive meshing with element size $\le 0.2\lambda_0$ (and $\le 0.05\lambda_0$ nearby the impedance surface) is utilized, yielding $\sim 100$ million degrees of freedom, and the problem solution is carried out via the MUMPS solver (with default parameters).

\bibliographystyle{IEEEtran}
\bibliography{SD_MoM_LWbibliography}%

\end{document}